\documentclass[journal]{IEEEtran}
\newcommand{\beq}{\begin{equation}}
\newcommand{\eeq}{\end{equation}}
\newcommand{\bsq}{\begin{subequations}}
\newcommand{\esq}{\end{subequations}}
\newcommand{\bq}{\begin{eqnarray}}
\newcommand{\eq}{\end{eqnarray}}
\newcommand{\bqn}{\begin{eqnarray*}}
\newcommand{\eqn}{\end{eqnarray*}}
\usepackage{cite}
\usepackage{url}
\usepackage{empheq}
\usepackage{float}
\usepackage{dsfont}
\usepackage[hidelinks]{hyperref}

\usepackage{algorithm}
\usepackage{algorithmic}

\usepackage{amsmath,amssymb,amsfonts}
\allowdisplaybreaks[4]
\usepackage{graphicx}
\usepackage{textcomp}
\usepackage{amsmath}
\usepackage{bbm}
\usepackage{enumerate}
\usepackage{epstopdf}
\usepackage{array}
\usepackage{booktabs}
\usepackage{subfigure}
\usepackage{multirow}
\usepackage{soul, color}
\usepackage[usenames,dvipsnames]{xcolor}
\usepackage[latin1]{inputenc}

\newtheorem{theorem}{Theorem}
\newtheorem{lemma}{Lemma}

\newtheorem{definition}{Definition}

\soulregister\cite7 
\soulregister\citep7 
\soulregister\citet7 
\soulregister\ref7
\soulregister\pageref7

\begin{document}

%
\title{Robust Generation Dispatch with Purchase of Renewable Power and Load Predictions}
%
%

\author{Rui Xie,
       Pierre Pinson,~\IEEEmembership{Fellow,~IEEE,}
       Yin Xu,~\IEEEmembership{Senior Member,~IEEE,}
       Yue Chen,~\IEEEmembership{Member,~IEEE}
\thanks{R. Xie and Y. Chen are with the Department of Mechanical and Automation Engineering, the Chinese University of Hong Kong, Hong Kong SAR. (email: ruixie@cuhk.edu.hk; yuechen@mae.cuhk.edu.hk)}
\thanks{P. Pinson is with the Dyson School of Design Engineering, Imperial College London, UK. (email: p.pinson@imperial.ac.uk)}
\thanks{Y. Xu is with the School of Electrical
Engineering, Beijing Jiaotong University, Beijing 100044, China. (email: xuyin@bjtu.edu.cn)}
}

%
%

\markboth{Journal of \LaTeX\ Class Files,~Vol.~XX, No.~X, Feb.~2019}%
{Shell \MakeLowercase{\textit{et al.}}: Bare Demo of IEEEtran.cls for IEEE Journals}
%



\maketitle

\begin{abstract}
The increasing use of renewable energy sources (RESs) and responsive loads has made power systems more uncertain. Meanwhile, thanks to the development of advanced metering and forecasting technologies, predictions by RESs and load owners are now attainable. Many recent studies have revealed that pooling the predictions from RESs and loads can help the operators predict more accurately and make better dispatch decisions. However, how the prediction purchase decisions are made during the dispatch processes needs further investigation. This paper fills the research gap by proposing a novel robust generation dispatch model considering the purchase and use of predictions from RESs and loads. The prediction purchase decisions are made in the first stage, which influence the accuracy of predictions from RESs and loads, and further the uncertainty set and the worst-case second-stage dispatch performance. This two-stage procedure is essentially a robust optimization problem with decision-dependent uncertainty (DDU). A mapping-based column-and-constraint generation (C\&CG) algorithm is developed to overcome the potential failures of traditional solution methods in detecting feasibility, guaranteeing convergence, and reaching optimal strategies under DDU. Case studies demonstrate the effectiveness, necessity, and scalability of the proposed model and algorithm.
\end{abstract}

\begin{IEEEkeywords}
robust generation dispatch, prediction purchase, decision-dependent uncertainty, mapping-based C\&CG
\end{IEEEkeywords}

%
\IEEEpeerreviewmaketitle

\section*{Nomenclature}
\addcontentsline{toc}{section}{Nomenclature}
\subsection{Abbreviations}
\begin{IEEEdescription}[\IEEEusemathlabelsep\IEEEsetlabelwidth{ssssssssssss}]
\item[C\&CG] Column-and-constraint generation
\item[DDU] Decision-dependent uncertainty
\item[DIU] Decision-independent uncertainty
\item[RES] Renewable energy resource
\item[RGD] Robust generation dispatch
\item[RO] Robust optimization
\end{IEEEdescription}

\subsection{Indices and Sets}
\begin{IEEEdescription}[\IEEEusemathlabelsep\IEEEsetlabelwidth{ssssssssssss}]
\item[$i \in \mathcal{I}_r, \mathcal{I}_d, \mathcal{I}$] Set of RESs/loads/agents
\item[$j \in \mathcal{J}$] Set of controllable generators
\item[$l \in \mathcal{L}$] Set of transmission lines
\item[$t \in \mathcal{T}$] Set of periods
\item[$\mathcal{X}$] Feasible set of the first-stage variable
\item[$\mathcal{X}_R$] Robust feasible set
\item[$\mathcal{U}$] Uncertainty set
\item[$\mathcal{Y}(x,u)$] Feasible set of the second-stage variable
\item[$V(\cdot)$] Vertex set
\end{IEEEdescription}
       
\subsection{Parameters}
\begin{IEEEdescription}[\IEEEusemathlabelsep\IEEEsetlabelwidth{ssssssssssss}]
\item[$I_r,I_d,I$] Number of RESs/loads/agents
\item[$J$] Number of generators.
\item[$L$] Number of transmission lines
\item[$T$] Number of periods
\item[$\overline{u}_{it}^r,\overline{u}_{it}^d$] Expected value of the maximum power output of RES $i$/power demand of load $i$ in period $t$
\item[$\sigma_{U_i}^2$] Variance of operator's estimation for the uncertainty of agent $i$
\item[$\delta,\xi$] Probability parameters of uncertainty set
\item[$m$] Prediction cost parameter
\item[$\theta_{jt}$] Binary parameter representing whether generator $j$ is on in period $t$
\item[$\rho_j$] Output cost coefficient of generator $j$
\item[$\rho_j^+,\rho_j^-$] Upward/downward regulation cost coefficient of generator $j$
\item[$\rho^c$] Penalty coefficient of real-time RES power curtailment
\item[$\gamma_j^+,\gamma_j^-$] Upward/downward reserve cost coefficient of generator $j$
\item[$R_j^+,R_j^-$] Maximum upward/downward reserve of generator $j$
\item[$\mathcal{R}_j^+,\mathcal{R}_j^-$] Maximum upward/downward ramping of generator $j$
\item[$P_j^{min}$, $P_j^{max}$] Minimum/maximum output of generator $j$
\item[$F_l$] Capacity of transmission line $l$
\item[$\pi_{jl},\pi_{il}$] Power transfer distribution factors
\end{IEEEdescription}

\subsection{Decision Variables}
\begin{IEEEdescription}[\IEEEusemathlabelsep\IEEEsetlabelwidth{ssssssssssss}]
\item[$u_{it}$] Uncertain power of agent $i$ in period $t$
\item[$\tau_i$] Prediction accuracy of the uncertainty of agent $i$
\item[$C_i$] Operator's payment to agent $i$ for information
\item[$p_{jt}$] Contemporary output of generator $j$ in $t$
\item[$r_{jt}^+,r_{it}^-$] Upward/downward reserve of generator $j$ in period $t$
\item[$p_{jt}^+,p_{jt}^-$] Upward/downward power adjustment of generator $j$ in period $t$
\item[$p_{it}^c$] Real-time power curtailment of RES $i$ in period $t$
\end{IEEEdescription}

\section{Introduction}

\IEEEPARstart{G}{eographically} distributed renewable energy sources (RESs) and responsive loads have flourished in recent years, posing great challenges on power system operations including higher risks of power imbalance and inadequate ramping capacities \cite{yang2021chance}. Robust generation dispatch (RGD) is an essential way to cope with the rising uncertainty \cite{wei2014robust}. It allows the transmission grid operator to minimize the operation cost under the worst-case uncertainty realizations. The performance of RGD largely depends on the quality of the uncertainty set it employs. A vast literature has been devoted to building more accurate uncertainty sets \cite{guan2013uncertainty}.


The parameters of an uncertainty set were determined to achieve the best trade-off between security and conservativeness \cite{bertsimas2004price}. A method for building polyhedral uncertainty sets based on the theory of coherent risk measures was introduced \cite{bertsimas2009constructing}. In addition to a better selection of parameters of the uncertainty set, enhancing the accuracy of uncertainty predictions is another important way.
An improved wind forecasting framework considering the spatio-temporal correlation in wind speed was developed and used to build a more accurate uncertainty set \cite{xie2013short}. A dynamic uncertainty set was proposed in \cite{lorca2014adaptive}, which improves the forecasts for the current period based on uncertainty realizations in the past few periods.
In the studies above, the uncertainty sets were estimated by the operator through processing their own data and forecasts. Thanks to the development of advanced metering and forecasting technologies, there is an emerging opportunity to further improve the forecast and uncertainty set. That is, through pooling of predictions from RES and load owners \cite{raja2023market}.


To aggregate and make use of the local predictions, data/information markets have been introduced. For example, regression markets were developed to aggregate local data for energy forecasting with proper incentives based on cooperative game theory \cite{pinson2022regression,han2022trading}. A data-sharing mechanism was designed for electricity retailers to improve their profits in the wholesale market \cite{wang2022mechanism}. A blockchain-based data transmission framework was developed for energy imbalance market \cite{cui2021blockchain}. Due to the high communication burdens and the risk of private data leakage, sharing data to perform a central prediction may not always be a good way. An alternative approach is to build an information market for aggregating predictions \cite{raja2023market,chen2020decentralized}. A binary prediction market was proposed in \cite{shamsi2021prediction} to forecast RESs. The day-ahead trading and valuation of load forecast were studied in \cite{sun2023trading}. This paper chooses to focus on the latter approach, i.e., to help the operator improve the uncertainty set by purchasing predictions from RESs and loads.

The fruitful works above focused on the forecasting tasks. But in fact, the forecast improvement via prediction purchase and the power system dispatch have mutual influences. Therefore, a holistic model that integrates them is necessary, which remains to be explored. This paper takes an initial step by proposing a robust generation dispatch model that allows the operator to purchase and use predictions from RESs and loads.
The proposed model turns out to be a case of robust optimization (RO) with decision-dependent uncertainty (DDU). The decision-dependent feature makes the traditional algorithms, such as Benders decomposition \cite{bertsimas2012adaptive} and column-and-constraint generation (C\&CG) \cite{zeng2013solving}, inapplicable. A reformulation method was proposed to solve RO with DDU \cite{nohadani2018optimization} and was extended to a more general uncertainty set \cite{lappas2018robust}. These two studies focused on static robust models, and the adaptive two-stage robust model is even more challenging. Modified Benders decomposition \cite{zhang2021robust,zhang2022nash}, adaptive C\&CG \cite{chen2022robust}, and multi-parametric programming \cite{avraamidou2020adjustable} methods were established to provide an exact solution. However, they might be time-consuming with an increasing number of resources. An efficient solution algorithm that is suitable for the DDU set studied in this paper is needed.

Our core objective is to provide a holistic model that allows the operator to exploit local predictions to build an improved uncertainty set and considers the tradeoff between prediction purchase cost and operation cost in the robust generation dispatch.
A mapping-based C\&CG algorithm is developed to solve the problem. Our contributions are two-fold:

(1) \emph{Robust Generation Dispatch Model with Purchase of Local Predictions}. A novel RGD model is developed to help the operator make better dispatch decisions by exploiting predictions purchased from agents (RESs and loads). Distinct from previous research that estimates the uncertainty set by the operator's own data/forecast, this paper builds an improved uncertainty set based on conditional expectations and variances derived by combining forecasts from the operator and agents. The proposed model turns out to be a RO with DDU. It is worth noting that the proposed model is generic and not limited to the RGD problem in this paper.



(2) \emph{Solution Algorithm}. A mapping-based C\&CG algorithm is developed to solve the proposed model by returning mapping constraints rather than the worst-case scenarios directly.  The proposed algorithm is proven to converge to the optimal solution within a finite number of iterations, while the traditional RO algorithms may fail to detect feasibility and guarantee optimality when dealing with DDU. Several comparisons are conducted to demonstrate its effectiveness.

The rest of this paper is organized as follows. Section \ref{sec-2} builds an improved uncertainty set based on the purchased predictions. The robust generation dispatch model and its solution algorithm are developed in Section \ref{sec-III} and \ref{sec:algorithm}, respectively. Numerical experiments are carried out in Section \ref{sec:case study}. Section \ref{sec:conclusion} concludes the paper.

\section{Improved Uncertainty Set with \\Purchase of Predictions}
\label{sec-2}
We start with a general form of two-stage RO and develop approaches to improving its uncertainty set via purchased predictions from agents. Later in Section \ref{sec-III}, the robust generation dispatch problem will be introduced in detail. The two-stage RO can be generally formulated as
\begin{align}\label{eq:tr-robust}
    \min_x ~ & \, \left\{f(x) + \max_{u \in \mathcal{U}} ~ \min_{y \in \mathcal{Y}(x,u)} g(y)\right\}, \\
    \mbox{s.t.}~ & \, x \in \mathcal{X} \cap \mathcal{X}_R, \nonumber
\end{align}
with
\begin{align}\label{eq:robust-feasible}
    \mathcal{X}_R=\{x~|~\mathcal{Y}(x,u) \ne \emptyset, \, \forall u \in \mathcal{U}\},
\end{align}
where $x$ and $y$ are the first- and second-stage decision variables, respectively. In parallel, $u$ represents the potential uncertainty realization varying within an uncertainty set $\mathcal{U}$. $\mathcal{X}$ and $\mathcal{X}_R$ are the feasible and robust feasible sets of $x$, respectively. According to \eqref{eq:robust-feasible}, a first-stage decision $x$ is robust feasible if and only if for any realization of uncertainty $u$, we can find a feasible second-stage solution $y \in \mathcal{Y}(x,u)$.

\begin{figure}[ht]
\centering
\includegraphics[width=0.95\columnwidth]{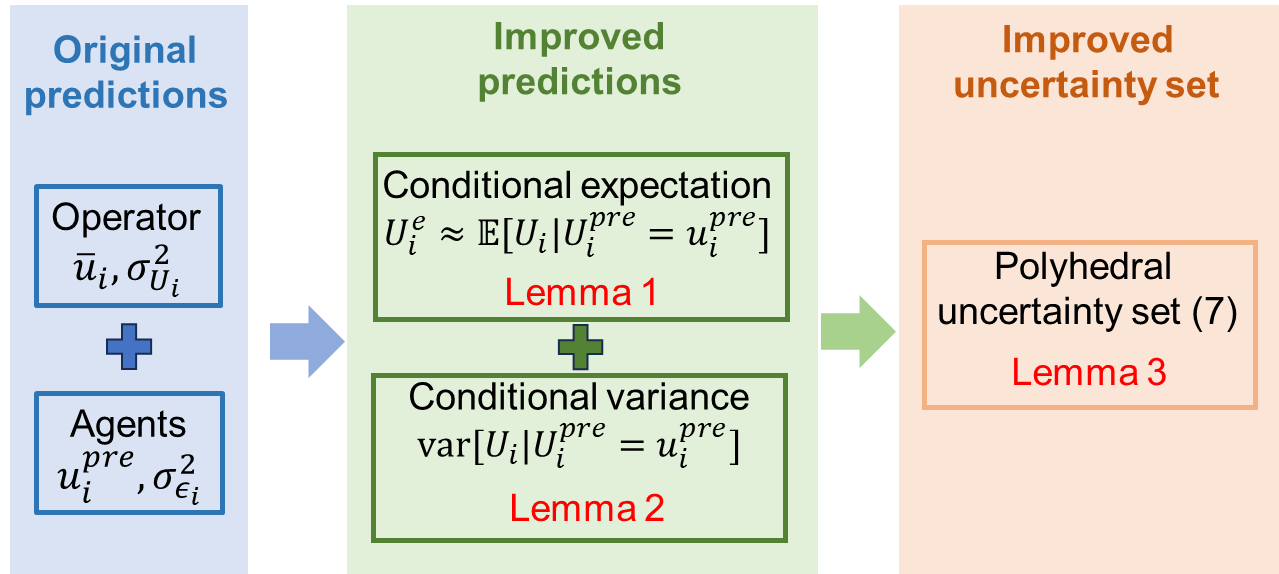}
\caption{Overview of Section \ref{sec-2}.}
\label{fig:procedure}
\end{figure}

Traditionally, the uncertainty set $\mathcal{U}$ is estimated by the operator based on their own forecasts. With the recent advance in metering and data analytic technologies, the agents (RESs and loads) can also produce their predictions \cite{shamsi2021prediction,sun2023trading}, which may be used to improve the operator's forecast. 
In the following, we first investigate how, if provided with predictions from agents, the operator can improve their forecasts and obtain an improved uncertainty set. To make it easier to follow, the procedures are summarized in Fig. \ref{fig:procedure}.

\subsection{Improved Forecasts} 
Suppose there are $I$ agents indexed by $i \in \mathcal{I}=\{1,2,...,I\}$. The uncertainty of agent $i \in \mathcal{I}$ is represented as a random variable $U_i$ in $\mathbb{R}$ with an unknown distribution, and $u_i$ is its realization. In this paper, we focus on the polyhedral uncertainty set \cite{wei2014robust} that is commonly used in RO. The two key parameters to determine such a set are the expectation and variance of the uncertain factor.
Traditionally, the operator estimates the expectation and variance of $U_i$ using their own historical data, denoted by $\overline{u}_i:=\mathbbm{E}[U_i]$ and $\sigma_{U_i}^2:=\mbox{var}(U_i)$, respectively. These two estimates may not be accurate enough due to the limited data of the operator. To predict more precisely, the operator can buy predictions from the agents. Suppose agent $i$'s forecast of $U_i$ is $U_i^{pre}$, then we have $U_i=U_i^{pre}+\epsilon_i$, where $\epsilon_i$ is a random noise.
The prediction $U_i^{pre}$ is also a random variable and let $u_i^{pre}$ be its realization. Denote by $\sigma_{\epsilon_i}^2:=\mbox{var}(\epsilon_i)$ the variance of $\epsilon_i$. The higher the agent's prediction accuracy, the smaller the $\sigma_{\epsilon_i}^2$. 
Agents are heterogeneous and may have different prediction errors reflected in $\sigma_{\epsilon_i}^2$.
Throughout the paper, we assume that: 

\vspace{0.5em}

\noindent \textbf{A1}:  $\{\epsilon_i,\forall i \in \mathcal{I}\}$ are independent. Each $\epsilon_i$ has zero expectation, i.e., $\mathbbm{E}[\epsilon_i]=0$, and $\epsilon_i$ is independent of $U_i$. 

\vspace{0.5em}

First, let us see how the operator can improve their estimation of $U_i$ with the help of prediction $u_i^{pre}$. 
We propose to use the conditional probability $\mathcal{P}(U_i|U_i^{pre}=u_i^{pre})$ as an approximation of the actual probability $\mathcal{P}(U_i)$ of $U_i$. Then the uncertainty set can be constructed based on the conditional expectation $\mathbbm{E}[U_i|U_i^{pre}=u_i^{pre}]$ and variance $\mbox{var}(U_i|U_i^{pre}=u_i^{pre})$. Generally, these two parameters can be complicated nonlinear functions of $u_i^{pre}$. For simplicity, in this paper, we adopt the \emph{best linear predictor} of $U_i$, $U_i^e := \alpha_i + \beta_i U_i^{pre}$, that minimizes the squared error expectation as follows, 
\begin{align}\label{eq:linear-predictor}
    \min_{\alpha_i, \beta_i}~ \mathbbm{E}\left[ \left(U_i - (\alpha_i + \beta_i U_i^{pre}) \right)^2 \right],
\end{align}
where $\alpha_i$ and $\beta_i$ are parameters to be determined. Denote the realization of $U_i^e$ by $u_i^e$.

\begin{lemma}\label{lemma-1} When A1 holds, the two parameters of the best linear predictor $U_i^e = \alpha_i + \beta_i U_i^{pre}$ are
\begin{align}\label{eq:expectation}
    \beta_i=\frac{\sigma^2_{U_i}}{\sigma^2_{U_i}+\sigma^2_{\epsilon_i}}, ~\alpha_i=(1-\beta_i)\overline{u}_i.
\end{align}
Moreover, $\mathbbm{E}[U_i-U_i^e] = 0$ and $\mathrm{cov}(U_i-U_i^e,U_i^{pre})=0$.
\end{lemma}

The proof of Lemma \ref{lemma-1} can be found in Appendix \ref{apendix-A}. 
When the error is very small ($\sigma_{\epsilon_i}^2 \to 0$), we have $\alpha_i=0, \beta_i=1$, and thus $u_i^e = u_i^{pre}$. It means that the prediction $u_i^{pre}$ is accurate so the operator just uses it. When the error is very high ($\sigma_{\epsilon_i}^2 \to \infty$), we have $\alpha_i=\bar{u}_i, \beta_i=0$, and thus $u_i^e=\overline{u}_i$. It means that the prediction $u_i^{pre}$ is so inaccurate that the operator cannot get a better estimation than the original one $\overline{u}_i$. Let $\eta_i:=U_i-U_i^e$. By Lemma \ref{lemma-1}, we have $\mathrm{cov}(\eta_i,U_i^{pre})=0$. So it is reasonable to make a stronger assumption:

\vspace{0.5em}

\noindent \textbf{A2}: $\eta_i$ and $U_i^{pre}$ are independent.

\vspace{0.5em}


Based on the best linear predictor, the conditional variance $\mbox{var}(U_i|U_i^{pre}=u_i^{pre})$ can be calculated by Lemma \ref{lemma-2}. 

\begin{lemma}\label{lemma-2} When A1 and A2 hold, we have
\begin{align}
    \mbox{var}(U_i|U_i^{pre}=u_i^{pre})=(1-\beta_i)^2\sigma^2_{U_i}+\beta_i^2\sigma_{\epsilon_i}^2
\end{align}
and $\mbox{var}(U_i|U_i^{pre}=u_i^{pre}) \le \sigma_{U_i}^2$.
\end{lemma}

The proof of Lemma \ref{lemma-2} can be found in Appendix \ref{apendix-B}. From Lemma \ref{lemma-1}, we have observed that the more accurate the prediction (i.e., the smaller the $\sigma^2_{\epsilon_i}$), the larger the $\beta_i$ (weight on $u_i^{pre}$). In an extreme case when $\sigma^2_{\epsilon_i}$ is zero, we have $u_i^e=u_i^{pre}$ and $\mbox{var}(U_i|U_i^{pre}=u_i^{pre})=0$. This indicates that with the prediction from agent $i$, the operator can know the exact value of $u_i$, and so there is no uncertainty. On the contrary, when $\sigma^2_{\epsilon_i} \to \infty$, we have $u_i^e=\overline{u}_i$, and $\mbox{var}(U_i|U_i^{pre}=u_i^{pre})=\sigma^2_{U_i}$, meaning that the agent's prediction is so inaccurate that the operator still uses the same estimation as if there were no purchased prediction. From the analysis above, we find that $\beta_i \in [0,1]$ can be a good indicator of prediction accuracy. We borrow similar concepts from economics and define the prediction accuracy as follows.

\begin{definition} (Prediction Accuracy \cite{vives2010information}) The parameter $\beta_i$ in \eqref{eq:expectation} can be formally defined as the prediction accuracy. To differentiate, we use $\tau_i$ to denote the prediction accuracy, i.e.,
\begin{align}\label{eq:accuracy}
    \tau_i:=\frac{\sigma^2_{U_i}}{\sigma^2_{U_i}+\sigma^2_{\epsilon_i}} \in [0,1].
\end{align}
\end{definition}

The accuracy of each agent's prediction to the operator is influenced both by the agent's forecasting technology and the incentive paid by the operator (will be explained later).

\subsection{Improved Uncertainty Set}

With the improved forecasts above, the operator can then construct an improved uncertainty set. 
Suppose there are $T$ periods indexed by $t \in \mathcal{T}=\{1,...,T\}$, then we have the best linear predictor $u_{it}^e$ and variance $\mbox{var}(U_{it}|U_{it}^{pre}=u_{it}^{pre})$ for all $i \in \mathcal{I}$ and $t \in \mathcal{T}$. We adopt a polyhedral uncertainty set with the following form.
\begin{align} \label{eq:uncertainty-set}
    \mathcal{U}(\tau)=~ &\{u_{it},\forall i \in \mathcal{I},\forall t \in \mathcal{T}: \nonumber\\
    ~ & u_{it}^0-u_{it}^h \le u_{it} \le u_{it}^0+u_{it}^h,\forall i \in \mathcal{I},\forall t \in \mathcal{T} \nonumber\\
    ~ & \sum \nolimits_i \frac{|u_{it}-u_{it}^0|}{u_{it}^h}\le \Gamma_S,~\forall t \in \mathcal{T}\nonumber\\
    ~ & \sum \nolimits_t \frac{|u_{it}-u_{it}^0|}{u_{it}^h}\le \Gamma_T,~\forall i \in \mathcal{I}\},
\end{align}
where $\Gamma_S$ and $\Gamma_T$ are the uncertainty budgets to restrain the spatial and temporal deviations from the forecast $u_{it}^0$. $u_{it}^0, u_{it}^h,\forall i,\forall t$ and $\Gamma_S,\Gamma_T$ are parameters to be determined based on the best linear predictor $u_{it}^e$ and variance $\mbox{var}(U_{it}|U_{it}^{pre}=u_{it}^{pre})$. Their values are given in Lemma \ref{lemma-3}.
Denote $v_{it}:=|u_{it}-u_{it}^0|/u_{it}^h,\forall i,\forall t$, so $v_{it}=|\eta_{it}/u_{it}^h|$ when $u_{it}^0 = u_{it}^e$. For simplicity, we assume that:

\vspace{0.5em}

\noindent \textbf{A3}:  $\{\eta_{it}/\sqrt{\mbox{var}[\eta_{it}]},\forall i,\forall t\}$ are independent and identically distributed (i.i.d.).

\vspace{0.5em}




\begin{lemma}\label{lemma-3} When A1-A3 hold, if the parameters of the uncertainty set \eqref{eq:uncertainty-set} are chosen as
\bsq
\begin{align}
    \label{eq:uncertainty-set-1}
    u_{it}^0=~ & u_{it}^e=(1-\tau_{it}) \overline u_{it} + \tau_{it} u_{it}^{pre},~\forall i,\forall t, \\
    u_{it}^h=~ & \sqrt{\mbox{var}[U_{it}|U_{it}^{pre}=u_{it}^{pre}]/(1-\delta)} = \sqrt{\mbox{var}[\eta_{it}]/(1-\delta)} \nonumber\\
    \label{eq:uncertainty-set-2}
    = ~ & \sqrt{\left[(1-\tau_{it})^2 \sigma_{U_i}^2 + \tau_{it}^2 \sigma_{\epsilon_i}^2\right]/(1-\delta)},~ \forall i,\forall t,\\
    \label{eq:uncertainty-set-3}
    \Gamma_S = ~ & \sqrt{\frac{I(1-\delta)(1+I-I \xi)}{1-\xi}},   \\
    \label{eq:uncertainty-set-4}
    \Gamma_T = ~ &  \sqrt{\frac{T(1-\delta)(1+T-T \xi)}{1-\xi}}.
\end{align}
\esq
Then, we can ensure that $\mathcal{P}(v_{it}\ge 1)\le 1-\delta, \forall i, \forall t$,  $\mathcal{P}\left(\sum_i v_{it} \ge \Gamma_S\right)\le 1-\xi, \forall t$ and $\mathcal{P}\left(\sum_t v_{it} \ge \Gamma_T\right) \le 1-\xi,\forall i$.
\end{lemma}

The proof of Lemma \ref{lemma-3} can be found in Appendix \ref{apendix-C}. 
With the help of purchased predictions, a smaller and more accurate uncertainty set can be obtained.
The improved uncertainty set \eqref{eq:uncertainty-set} depends on the prediction accuracy $\tau_{it}, \forall t$ which is influenced by the payment of the operator to agent $i$ for buying the agent's prediction (denoted by $C_i,\forall i$). The higher the payments, the more accurate the predictions. The relationship between the operator's payment and the prediction accuracy is given in the next section. The improved uncertainty set is impacted by the payment $C_i,\forall i$ determined in the first stage, and thus, is decision-dependent.

\section{Robust Generation Dispatch Model With Purchase of Predictions}
\label{sec-III}
With the improved uncertainty set above, in the following, we develop the RGD model integrating the purchase and use of predictions from agents (RESs and loads).

\subsection{General Form}
To integrate the prediction purchase and RGD processes, we need to derive the relationship between the operator's payment $C_i$ and the prediction accuracy $\tau_{it}$. To this end, we first quantify the prediction cost for agent $i$ as
\begin{align}
    h_i(\tau_{i})=\sum_{t=1}^T \frac{m}{\sigma_{\epsilon_i}^2} =\sum_{t=1}^T \frac{ m}{\sigma_{U_{it}}^2}\frac{\tau_{it}}{1-\tau_{it}},
\end{align}
where $m$ is a given cost parameter. The more accurate the prediction (the smaller the $\sigma_{\epsilon_i}^2$), the higher the cost $h_i(\tau_i)$. For notation conciseness, we assume that $\tau_{i1}=\tau_{i2}=...=\tau_{iT}=:\tau_i$ and $\sigma_{U_{i1}}^2=\sigma_{U_{i2}}^2=...=\sigma_{U_{iT}}^2=:\sigma_{U_i}^2$. But it is worth noting that the proposed model and algorithm can also be applied to the cases with heterogeneous $\tau_{it},\forall t$ and $\sigma_{U_{it}}^2,\forall t$. Let $\hat m=T m$, then
\begin{align}\label{eq:prediction-cost}
    h_i(\tau_i)=\frac{\hat m}{\sigma_{U_i}^2}\frac{\tau_i}{1-\tau_i}.
\end{align}
The function \eqref{eq:prediction-cost} showing how the cost of prediction changes with the accuracy is commonly used in economics \cite{vives2010information}.
To ensure that the agents are willing to provide their predictions, the payment should be able to cover their prediction cost, i.e.,
\begin{align}\label{eq:precost}
    C_i \ge h_i(\tau_i)=\frac{\hat m}{\sigma_{U_i}^2}\frac{\tau_i}{1-\tau_i},\forall i \in \mathcal{I}.
\end{align}
\begin{definition} (Value of Prediction)
The operator's payment $C_i$ to agent $i$ can be formally defined as the value of prediction from agent $i$.
\end{definition}

The total payment for buying predictions from agents is also a cost of the operator in the first stage. Therefore, the two-stage RO model considering the purchase and use of predictions can be formulated as
\begin{align}\label{eq:new-robust}
     \min_{x,C,\tau} ~ & \left\{f(x) +\sum_{i \in \mathcal{I}} C_i + \max_{u \in \mathcal{U}(\tau)} ~ \min_{y \in \mathcal{Y}(x,u)} g(y)\right\}, \\
    \mbox{s.t.}~ & x \in \mathcal{X} \cap \tilde{\mathcal{X}}_R, \nonumber\\
    ~ & \eqref{eq:precost}, ~  0 \le \tau_i \le 1,\forall i \in \mathcal{I}, \nonumber
\end{align}
where
\begin{align}\label{eq:robust-feasible-new}
    \tilde{\mathcal{X}}_R=\{x~|~\mathcal{Y}(x,u) \ne \emptyset,\forall u \in \mathcal{U}(\tau)\}.
\end{align}

As the focus of this paper is the operator's decision-making in RGD, the predictive information market design between the operator and agents will be left for future study. 
The model \eqref{eq:new-robust} is an RO with DDU since the uncertainty set $\mathcal{U}(\tau)$ is influenced by the first-stage decision $\tau$. The details of the objective function and constraints are given below.


\subsection{Detailed Robust Generation Dispatch Model}
\label{sec:dispatch_example}
There are $J$ controllable generators indexed by $j \in \mathcal{J}=\{1,...,J\}$, $L$ lines indexed by $l \in \mathcal{L} = \{1,...,L\}$, $I_r$ RESs indexed by $i \in \mathcal{I}_r=\{1,...,I_r\}$, and $I_d$ loads indexed by $i \in \mathcal{I}_d=\{I_r+1,...,I_r+I_d\}$ in a transmission grid. Let $\mathcal{I} = \mathcal{I}_r \cup \mathcal{I}_d$ denote the set of agents whose power outputs/demands $u_{it},\forall i \in \mathcal{I},\forall t\in \mathcal{T}$ are uncertain. In particular, for $i \in \mathcal{I}_r$, $u_{it}$ represents the uncertain maximum power output of RES $i$ in period $t$; for $i \in \mathcal{I}_d$, $u_{it}$ represents the uncertain power demand of load $i$ in period $t$. The operator may buy predictions from the RESs and loads to improve their predictions. The best linear predictors $u_{it}^e$ of $u_{it}$ can be obtained by \eqref{eq:uncertainty-set-1}. 

In the first stage (day-ahead pre-dispatch stage), the transmission grid operator decides on the reference output and reserve capacity of the controllable generators and the payments for buying predictions from RESs and loads. In the second stage (re-dispatch stage), knowing the exact RES power outputs and exact demands, the operator adjusts the output of controllable generators within their reserve capacity or curtails RES power to maintain power balance. The robust generation dispatch (RGD) problem can be formulated as
\bsq \label{eq:RS-UL}
\begin{align}
   & \min_{p,r^\pm, \atop C,\tau} ~  \underbrace{\sum_{t \in \mathcal{T}} \sum_{j \in \mathcal{J}} (\varrho_j p_{jt}+\gamma_j^+ r_{jt}^+ + \gamma_j^- r_{jt}^-) }_{f(x)}+ \sum_{i \in \mathcal{I}} C_i \nonumber\\
    & + \max_{u \in \mathcal{U}(\tau)} \min_{(p^\pm,p^c)\in \atop \mathcal{Y}(p,r^\pm,u)} \underbrace{\sum_{t \in \mathcal{T}} \left(\sum_{j \in \mathcal{J}} (\varrho_j^+ p_{jt}^+ + \varrho_j^- p_{jt}^-)+\sum_{i \in \mathcal{I}_r} \rho^c p_{it}^c\right)}_{g(y)}, \label{eq:RS-UL.1}\\
    & \mbox{s.t.}~  (p,r^\pm) \in \mathcal{X} \cap \tilde{\mathcal{X}}_R, \\
    ~ & C_i \ge \frac{m}{\sigma^2_{U_i}}\frac{\tau_i}{1-\tau_i},~0 \le \tau_i \le 1, ~\forall i \in \mathcal{I}, \label{eq:RS-UL.7} 
\end{align}
\esq
where
\bsq
\begin{align}
    & \mathcal{X}= \Big\{ (p,r^\pm) ~| \nonumber\\
     & \sum \nolimits_{j \in \mathcal{J}} p_{jt}+\sum \nolimits_{i \in \mathcal{I}_r} u_{it}^e=\sum \nolimits_{i \in \mathcal{I}_d}  u_{it}^e,\forall t \in \mathcal{T}, \label{eq:RS-UL.2}\\
     & 0 \le r_{jt}^+ \le R_j^+ \theta_{jt}, 0 \le r_{jt}^- \le R_j^- \theta_{jt}, \forall j \in \mathcal{J},\forall t \in \mathcal{T}, \label{eq:RS-UL.3}\\
     & P_j^{min}\theta_{jt}+r_{jt}^- \le p_{jt} \le P_{j}^{max}\theta_{jt}-r_{jt}^+, \forall j \in \mathcal{J},\forall t \in \mathcal{T},
     \label{eq:RS-UL.4}\\
     & (p_{jt}+r_{jt}^+) -(p_{j(t-1)}-r_{j(t-1)}^-) \le \mathcal{R}_j^+\theta_{j(t-1)} \nonumber\\ &   + P_j^{max}(1-\theta_{j(t-1)}), \forall j \in \mathcal{J},\forall t=2,...,T, \label{eq:RS-UL.5}\\
     & - (p_{jt}-r_{jt}^-) + (p_{j(t-1)}+r_{j(t-1)}^+) \leq \mathcal{R}_j^- \theta_{jt} \nonumber\\
     & + P_j^{max}(1-\theta_{jt}),\forall j \in \mathcal{J},\forall t=2,...,T, \label{eq:RS-UL.6}\\
     \! & \!-F_l \!\le\! \sum_{j \in \mathcal{J}} \pi_{jl}p_{jt}\! + \! \sum_{i \in \mathcal{I}_r} \pi_{il}u_{it}^e\! -\! \sum_{i \in \mathcal{I}_d} \pi_{il}u_{it}^e \!\le\! F_l,\!\forall l,\!\forall t \Big\},  \label{eq:RS-UL.9}
\end{align}
\esq
and
\bsq \label{eq:RS-LL}
\begin{align}
    & \mathcal{Y}(p,r^\pm,u)= \Big\{ p^\pm,p^c~| \nonumber\\
    ~ & 0 \le p_{jt}^+ \le r_{jt}^+, 0 \le p_{jt}^- \le r_{jt}^-,\forall j \in \mathcal{J},\forall t \in \mathcal{T}, \label{eq:RS-LL.1}\\
    ~ & \sum_{j \in \mathcal{J}} (p_{jt}+ p_{jt}^+ - p_{jt}^-) + \sum_{i \in \mathcal{I}_r} (u_{it}-p_{it}^c) = \sum_{i \in \mathcal{I}_d} u_{it},\forall t \in \mathcal{T}, \label{eq:RS-LL.2}\\
    ~ & 0 \leq p_{it}^c \leq u_{it}, \forall i \in \mathcal{I}_r, \forall t \in \mathcal{T}, \label{eq:RS-LL.3}\\
     ~ & -F_l \le \sum_{j \in \mathcal{J}} \pi_{jl} (p_{jt}+ p_{jt}^+ - p_{jt}^-) + \sum_{i \in \mathcal{I}_r} \pi_{il}(u_{it}-p_{it}^c)  \nonumber\\
    ~ & - \sum_{i \in \mathcal{I}_d} \pi_{il}u_{it} \leq F_l,\forall l \in \mathcal{L},\forall t \in \mathcal{T}\Big\}.\label{eq:RS-LL.4}
\end{align}
\esq
The objective function \eqref{eq:RS-UL.1} minimizes the total cost under the worst-case scenario, i.e., the total generation-related cost and prediction purchase payment in the first stage plus the total generation adjustment cost and curtailment penalty in the second stage. In the first stage, the decision variable $x$ consists of the reference output $\{p_{jt},\forall j,\forall t\}$ and upward/downward reserve capacity $\{r_{jt}^+,r_{jt}^-,\forall j,\forall t\}$ of controllable generators. $C_i,\forall i$ is the payment for buying information from the agents. The second-stage decision variable $y$ includes the upward/downward power output adjustment $\{p_{jt}^+,p_{jt}^-,\forall j,\forall t\}$ of controllable units and the real-time RES power curtailment $\{p_{it}^c,\forall i \in \mathcal{I}_r,\forall t\}$. $\theta_{jt}$ is a binary parameter: $\theta_{jt}=1$ means in period $t$, generator $j$ is on and otherwise $\theta_{jt}=0$. $\varrho_j$, $\gamma_j^\pm$, and $\varrho_j^\pm$ are the cost coefficients of power output, upward/downward reserve, and upward/downward regulation, respectively. $\rho^c$ is the real-time curtailment penalty coefficient.

Constraints \eqref{eq:RS-UL.2}-\eqref{eq:RS-LL.4} stipulate the operational limits. Constraints \eqref{eq:RS-UL.2} and \eqref{eq:RS-LL.2} are the power balance conditions. The upward/downward reserve capacity should not exceed the bounds $R_j^\pm$ as in \eqref{eq:RS-UL.3}. The upper/lower power limits of controllable generators considering reserve requirements are given in \eqref{eq:RS-UL.4}, where $P_j^{min}/P_j^{max}$ is the minimum/maximum power output. \eqref{eq:RS-UL.5}-\eqref{eq:RS-UL.6} ensure the satisfaction of ramping limits when offering reserves \cite{chen2022robust}. The upper/lower power limits of RES curtailment are given in \eqref{eq:RS-LL.3}. The network capacity limits are imposed in \eqref{eq:RS-UL.9} and \eqref{eq:RS-LL.4}; $F_l$ is the power flow limit of line $l$ and $\pi_{il}, \pi_{jl}$ are the power transfer distribution factors (PTDFs) deduced from the DC power flow model, so constraints \eqref{eq:RS-UL.2}, \eqref{eq:RS-UL.9}, \eqref{eq:RS-LL.2}, and \eqref{eq:RS-LL.4} constitute the network model. Constraint \eqref{eq:RS-LL.1} ensure that the power adjustment is within the reserve capacity.

\emph{Remark:} Load shedding is not allowed in the proposed model \eqref{eq:RS-UL}-\eqref{eq:RS-LL}. This is because in some countries such as China, load shedding is viewed as an operation failure and is prohibited \cite{wei2014robust}.  But it is worth noting that the proposed model can also accommodate load shedding by adding slack variables to \eqref{eq:RS-LL.2}, \eqref{eq:RS-LL.4} and a penalty term in the objective function.


As mentioned earlier, the proposed model \eqref{eq:RS-UL}-\eqref{eq:RS-LL} is an RO with DDU. The traditional algorithms such as Benders decomposition and C\&CG cannot be directly applied since they may fail to converge or lead to suboptimal solutions. In the next section, a mapping-based C\&CG algorithm will be developed to overcome this difficulty.


\section{Solution Algorithm}
\label{sec:algorithm}
In this section, a mapping-based C\&CG algorithm is developed to solve the problem \eqref{eq:RS-UL}, an RO with DDU. Notice that the re-dispatch problem \eqref{eq:RS-LL} is a linear program, $g(y)$ and $\mathcal{Y}(x,u)$ can be expressed by
\begin{align}
    g(y)=~ & c^{\top}y, \\
    \mathcal{Y}(x,u)=~ & \{y \in \mathbb{R}^{n_y}~|~ Ax+By+Du \le q\}.
\end{align}

\subsection{Second-Stage Problem Transformation}
Given the first-stage decision $x \in \mathcal{X}$ and $\tau_i \in [0,1],\forall i \in \mathcal{I}$, the second-stage problem is a bilevel optimization:
\begin{align}
\label{eq:second-stage}
    S(x,\tau) = \max_{u \in \mathcal{U}(\tau)} \min_{y \in \mathcal{Y}(x,u)} c^\top y,
\end{align}
which is equivalent to the sub-problem (\textbf{SP}) \eqref{eq:SP} by converting the inner ``min'' problem to its KKT condition.
\bsq\label{eq:SP}
\begin{align}
    \textbf{SP}:~ & \max_{u \in \mathcal{U}(\tau), y, \nu} c^{\top} y, \label{eq:SP-1}\\
    & \mbox{s.t.}~ B^{\top}\nu=c,  \label{eq:SP-2}\\
    ~ & 0 \le -\nu \perp \left[-(Ax+By+Du)+q\right] \ge 0. \label{eq:SP-3}
\end{align}
\esq
The complementary slackness condition \eqref{eq:SP-3} can be linearized by the Big-M method \cite{pineda2019solving}.

Furthermore, for a given first-stage decision $(x,\tau)$, the problem \eqref{eq:second-stage} may be infeasible. Remember that we need to ensure $x$ is robust feasible ($x \in \tilde{\mathcal{X}}_R$), so we construct the following relaxed problem for checking feasibility.
\bsq\label{eq:second-stage-relax}
\begin{align}
    F(x,\tau) = \max_{u \in \mathcal{U}(\tau)}~ &  \min_{y,s} 1^{\top}s, \\
    \mbox{s.t.} ~&  Ax+By+Du-s \le q, s\ge 0.
\end{align}
\esq

Obviously, the relaxed problem \eqref{eq:second-stage-relax} is always feasible. Moreover, we have the original problem \eqref{eq:second-stage} is feasible if and only if $s$ is an all-zero vector at the optimum. Similarly, the relaxed problem \eqref{eq:second-stage-relax} is equivalent to
\bsq\label{eq:feasibility-check}
\begin{align}
    \textbf{FC}:~& \max_{u \in \mathcal{U}(\tau) \atop y, s, \nu, \mu}~ 1^{\top} s, \\
    & \mbox{s.t.} ~ B^{\top}\nu =0, \\
    ~ & -\nu+\mu=1,\\
    ~ & 0 \le -\nu \perp \left[-(Ax+By+Du-s)+q\right] \ge 0,\\
    ~ & 0 \le \mu \perp s \ge 0.
\end{align}
\esq
We call \eqref{eq:feasibility-check} the feasibility-check (\textbf{FC}) problem.

Given a candidate first-stage decision, we first solve the \textbf{FC} problem to check whether $x \in \tilde{\mathcal{X}}_R$. If not, a feasibility cut will be returned; otherwise, we continue to solve the \textbf{SP} problem to identify an optimality cut.

\begin{lemma}\label{lemma-4}
Suppose $u^*$ is the optimal solution of \textbf{SP} or \textbf{FC}, then $u^*$ can be reached at a vertex of $\mathcal{U}(\tau)$. 
\end{lemma}

The proof of Lemma \ref{lemma-4} is similar to that in \cite{konno1976cutting} and is omitted here. The traditional RO algorithms return the worst-case scenario $\{u_{it}^*,\forall i,\forall t\}$ directly to the master problem to generate a feasibility/optimality cut. However, when dealing with DDU, a previously selected scenario may no longer be a vertex of the new uncertainty set when the first-stage decision changes ($\mathcal{U}(\tau)$ changes with $\tau$). This causes the traditional algorithms to fail to find the optimal solutions.

\begin{figure}[ht]
\centering
\includegraphics[width=0.5\columnwidth]{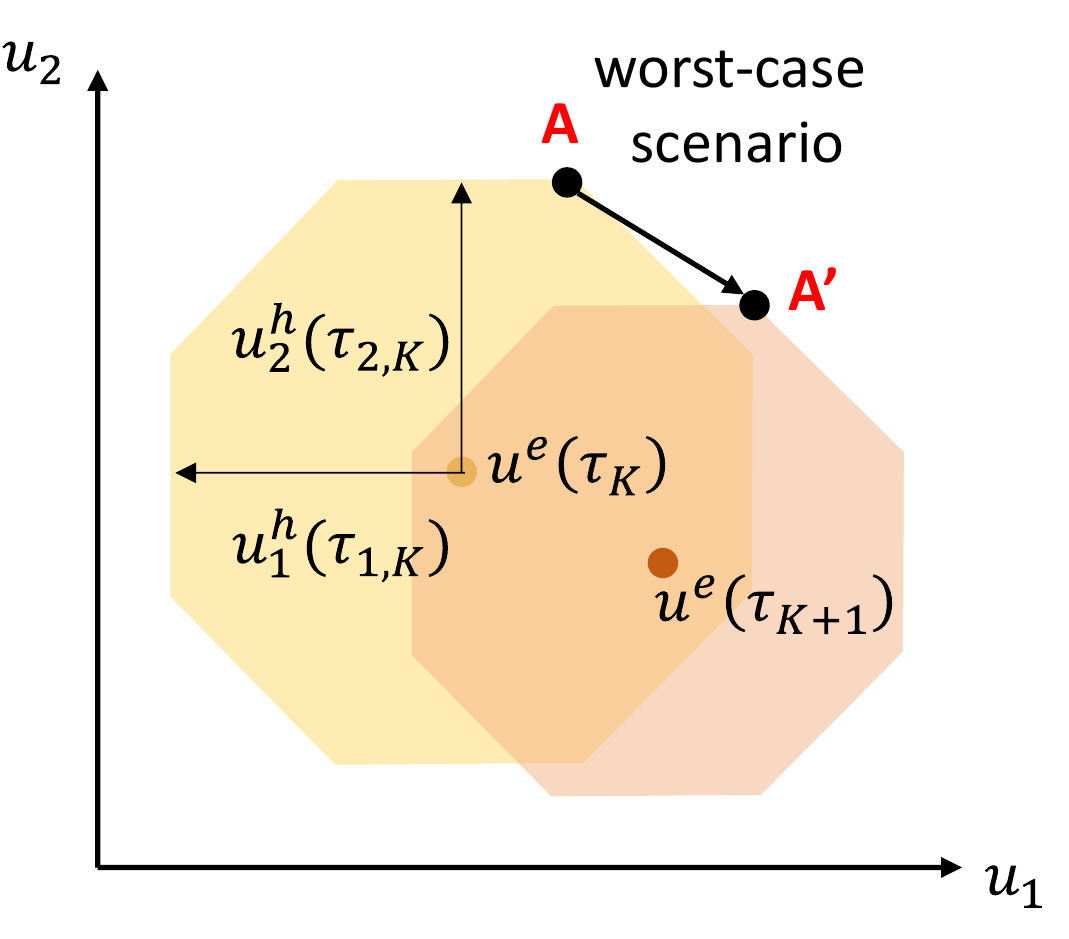}
\caption{Illustration of \eqref{eq:nonlinear-constraint}. The yellow and orange regions are the uncertainty sets in the $K$ and $K+1$ iterations, respectively. Point A is the worst-case scenario in the $K$ iteration. When $u^e$ and $u^h$ change with $\tau$, \eqref{eq:nonlinear-constraint} moves point A to point $\mbox{A}'$, which is a vertex of the new uncertainty set.}
\label{fig:algorithm}
\end{figure}

To tackle this problem, instead of returning the scenario $u_{it}^*,\forall i,\forall t$ directly, we propose to map the worst-case scenario to an unchanged vertex set and return the mapping constraints. To be specific, $\mathcal{U}(\tau)$ can be equivalently represented as
\begin{align}
    \mathcal{U}(\tau)=~ &\left\{u_{it},\forall i \in \mathcal{I},\forall t \in \mathcal{T}~| \right. \nonumber\\
    ~ &  u_{it}=u_{it}^e(\tau_i) + u_{i}^h(\tau_i) \phi_{it},\forall i \in \mathcal{I},\forall t \in \mathcal{T}, \nonumber\\
    ~ & \left.\{\phi_{it}, \forall i \in \mathcal{I}, \forall t \in \mathcal{T} \} \in \Phi \right\}, \nonumber
\end{align}
where
\begin{align}
    \Phi:= ~ &\{\phi_{it},\forall i \in \mathcal{I},\forall t \in \mathcal{T}~| \nonumber\\
    ~ & -1 \le \phi_{it} \le 1,\phi_{it} \le \psi_{it}, -\phi_{it} \le \psi_{it},\forall i \in \mathcal{I},\forall t \in \mathcal{T}, \nonumber\\
    ~ & \sum \nolimits_{i \in \mathcal{I}} \psi_{it} \le \Gamma_S, \forall t \in \mathcal{T}, \sum \nolimits_{t \in \mathcal{T}} \psi_{it} \le \Gamma_T, \forall i \in \mathcal{I} \}. \nonumber
\end{align}
Denote the vertex set of $\Phi$ by $V(\Phi)$, then the vertex set of $\mathcal{U}(\tau)$ can be represented as
\begin{align}
    V(\mathcal{U}(\tau))=~ &\{u_{it},\forall i \in \mathcal{I},\forall t \in \mathcal{T}~| \nonumber\\
    ~ &  u_{it}=u_{it}^e(\tau_i) + u_{i}^h(\tau_i) \phi_{it},\forall i \in \mathcal{I},\forall t \in \mathcal{T}, \nonumber\\
    ~ & \{\phi_{it}, \forall i \in \mathcal{I}, \forall t \in \mathcal{T} \} \in V(\Phi) \}, \nonumber
\end{align}
where $V(\Phi)$ is an unchanged vertex set. For a worst-case scenario $u^*$ generated by \textbf{FC} or \textbf{SP}, we can get the correspondingly $\phi_{it}^*,\forall i,\forall t$. Then, instead of returning $u^*$ to the master problem, we return the following mapping constraints:
\begin{align}
    \label{eq:nonlinear-constraint}
     u_{it}=u_{it}^e(\tau_i)+u_{i}^h(\tau_i)\phi_{it}^*,\forall i \in \mathcal{I},\forall t \in \mathcal{T}.
\end{align}
Here, 
$u_{it}$, $u_{it}^e(\tau_i)$, and $u_i^h(\tau_i)$ are all variables in the master problem. When $\tau_i$ changes, the $u_{it},\forall i,\forall t$ calculated by \eqref{eq:nonlinear-constraint} remains at a vertex of $\mathcal{U}(\tau)$, as illustrated in Fig. \ref{fig:algorithm}.

\subsection{Mapping-Based C\&CG Algorithm}
With the mapping constraints returned, the master problem (\textbf{MP}) can be formulated as
\bsq
\label{eq:MPv1}
\begin{align}
    \textbf{MP}: ~ & \min_{x, C,\tau, \zeta, y^k, u^k}~ f(x)+\sum_i C_i+\zeta, \label{eq:MPv1.1}\\
    & \mbox{s.t.}~ x \in \mathcal{X}, \eqref{eq:precost}, \tau_i \in [0,1],\forall i \in \mathcal{I}, \label{eq:MPv1.2}\\
    ~ & \zeta \ge c^{\top} y^k, \forall k \in [K],  \label{eq:MPv1.3}\\
    ~ & Ax+By^k+Du^k \le q, \forall k \in [K], \label{eq:MPv1.4} \\
    ~ & u^k=u^e(\tau) + u^h(\tau)\phi_k, \forall k \in [K], \label{eq:MPv1.5}
\end{align}
\esq
where the symbols with superscript $k$ are variables while the symbols with subscript $k$ are given parameters. $[K]$ represents all positive integers not exceeding $K$. 

\begin{algorithm}
\normalsize
\caption{{\bf Mapping-Based C\&CG Algorithm}}
\begin{algorithmic}[1]
\STATE \textbf{Initiation}: Error tolerance $\epsilon>0$; $K=0$; $UB_K = +\infty$.
 
\STATE \textbf{Solve the Master Problem}

Solve the $\textbf{MP}$ \eqref{eq:MPv1}. Derive the optimal solution $(x_{K+1}^*, C_{K+1}^{*}\tau_{K+1}^*, \zeta_{K+1}^*, y^{1*}, ..., y^{K*}, u^{1*},..., u^{K*})$ and update $LB_{K+1}=f(x_{K+1}^*)+\sum_{i \in \mathcal{I}} C_{i,K+1}^*+\zeta_{K+1}^*$.

\STATE \textbf{Solve the Feasibility-check Problem}

Solve the $\textbf{FC}$ \eqref{eq:feasibility-check} with $(x_{K+1}^{*},\tau_{K+1}^*)$.

Let $(u_{K+1}^*,\phi_{K+1}^*,\nu_{K+1}^*, \mu_{K+1}^*, y_{K+1}^*,s_{K+1}^*)$ be the optimal solution. If $1^{\top}s_{K+1}^*>0$, let $UB_{K+1} = UB_K$ and go to Step 6. Otherwise, go to Step 4.

\STATE \textbf{Solve the Sub-problem}

Solve the $\textbf{SP}$ \eqref{eq:SP} with $(x_{K+1}^{*},\tau_{K+1}^*)$. Denote the optimal solution by $(u_{K+1}^*,\phi_{K+1}^*,\nu_{K+1}^*, y_{K+1}^*)$. Let
\begin{align}
    UB_{K+1}=f(x_{K+1}^*)+\sum \nolimits_{i \in \mathcal{I}} C_{i,K+1}^*+c^{\top} y_{K+1}^* \nonumber
\end{align}

\STATE If $|UB_{K+1}-LB_{K+1}|\le \epsilon$, terminate and output $(x_{K+1}^{*},\tau_{K+1}^*)$. Otherwise, go to Step 6.

\STATE Create variables $(y^{K+1}, u^{K+1})$ and add the following constraints to \textbf{MP} \eqref{eq:MPv1}:
\begin{align}
    ~ & \zeta \ge c^{\top} y^{K+1},~  Ax+By^{K+1} + Cu^{K+1} \le q, \nonumber\\
    ~ & u^{K+1}=u^e(\tau)+u^h(\tau)\phi_{K+1}^*. \nonumber
\end{align}
Update $K=K+1$ and go to Step 2.

\end{algorithmic}
\end{algorithm} 

The overall procedure of the proposed mapping-based C\&CG algorithm is given in Algorithm 1. The proposed algorithm is different from the traditional C\&CG algorithm \cite{zeng2013solving} as it returns the mapping constraints instead of the worst-case scenarios $\{ u_{it}^*,\forall i, \forall t \}$ to the master problem to generate new cuts. To be specific, in the traditional C\&CG algorithm, Step 6 in Algorithm 1 is replaced by: ``6: Create variables $y^{K+1}$ and add the following constraints to \textbf{MP} \eqref{eq:MPv1}:
\begin{align}
	\zeta \ge c^{\top} y^{K+1}, ~Ax+By^{K+1} + Cu_{K+1}^* \le q. \nonumber
\end{align}
Update $K=K+1$ and go to Step 2." Note that $u_{K+1}^*$ is a constant vector obtained by solving \textbf{FC} \eqref{eq:feasibility-check} or \textbf{SP} \eqref{eq:SP}.

\begin{theorem}\label{thm-1}
Let $n_{U} := |V(\Phi)|$ be the number of extreme points of $\Phi$.
The mapping-based C\&CG algorithm generates the optimal solution to problem \eqref{eq:new-robust} within $\mathcal{O}(n_{U})$ iterations.
\end{theorem}

The proof of Theorem \ref{thm-1} can be found in Appendix \ref{apendix-D}. While the traditional algorithms fail to guarantee finite-step convergence and optimality of the obtained strategy, the proposed algorithm can overcome these limitations. Case studies in Section \ref{secV-C} also demonstrate the advantages and necessity of the proposed algorithm.

\subsection{Transformation and Linearization}

In Algorithm 1, the master problem \eqref{eq:MPv1} is highly nonlinear due to the term $u^e(\tau)$, $u^h(\tau)$, and the constraint \eqref{eq:precost}. In the following, we show how to turn \eqref{eq:MPv1} into a solvable form.

First, it is easy to prove that at the robust optimum, we have $C_i=h_i(\tau_i)$. Otherwise, if $C_i > h_i(\tau_i)$, we can always reduce $C_i$ a little bit without changing the value of the other variables, so that all constraints are still satisfied but the objective value decreases. This contradicts the definition of the robust optimum. Therefore, we can eliminate constraint \eqref{eq:precost} and replace $\sum_i C_i$ in the objective function with $\sum_i h_i(\tau_i)$.

Second, if we let $u_i^h,\forall i \in \mathcal{I}$ be the decision variables and use them to represent $\tau_i,\forall i \in \mathcal{I}$, then the prediction cost $\sum_i h_i(\tau_i)$ can be represented by
\begin{align}\label{eq:uh}
     \sum \nolimits_{i \in \mathcal{I}} h_i(\tau_i) =\sum \nolimits_{i \in \mathcal{I}} \hat m\left(\frac{1}{(1-\delta) (u_{i}^h)^2}-\frac{1}{\sigma_{U_i}^2}\right).
\end{align}
Let $\tilde{h}_i(u_i^h)$ denote the term in the right-hand side of \eqref{eq:uh} for each $i \in \mathcal{I}$, where $u_i^h$ is the decision variable.
Then, $\tilde{h}_i(u_i^h)$ is a convex function. Similarly, $u_{it}^e(\tau)$ can be represented by $u_i^h$, which is
\begin{align}\label{eq:ue}
    u_{it}^e(u_i^h)\! =\! \frac{(1-\delta)(u_{i}^h)^2}{\sigma_{U_i}^2} \bar{u}_{it} +  \left(\!1\!-\!\frac{(1-\delta)(u_{i}^h)^2}{\sigma_{U_i}^2}\!\right) u_{it}^{pre},\forall i,\forall t.
\end{align}

Next, we introduce a new variable $\tilde u_{it}^e,\forall i \in \mathcal{I},\forall t \in \mathcal{T}$, use it to replace $u^e(\tau)$ in \eqref{eq:MPv1.5}, and add the following penalty function to the objective:
\begin{align}
    \mathcal{H}(\tilde u^e, u^h) = \iota \sum \nolimits_{i \in \mathcal{I}} \sum \nolimits_{t \in \mathcal{T}} \left(\tilde u_{it}^e-u_{it}^e(u_{i}^h)\right)^2,
\end{align}
where $\iota$ is a large constant. Then, the remaining nonlinear term in the objective function, $\sum_i \tilde{h}_i(u_i^h) + \mathcal{H}(\tilde u^e, u^h)$, can be linearized by a convex combination approach \cite{wu2011tighter}, after which the \textbf{MP} problem has been turned into a linear program that can be solved efficiently.



\section{Case Studies}
\label{sec:case study}
We first use a simple 5-bus system to verify the proposed method and reveal some interesting phenomena; then, larger systems (33-, 69-, and 123-bus systems) are tested to show the scalability. Detailed data can be found in \cite{xie2022github}. The mapping-based C\&CG algorithm is implemented in MATLAB with GUROBI 9.5. All the simulations are conducted on a laptop with Intel i7-12700H processor and 16 GB RAM.

\subsection{Benchmark}
   
\begin{figure}[ht]
\centering
\includegraphics[width=0.9\columnwidth]{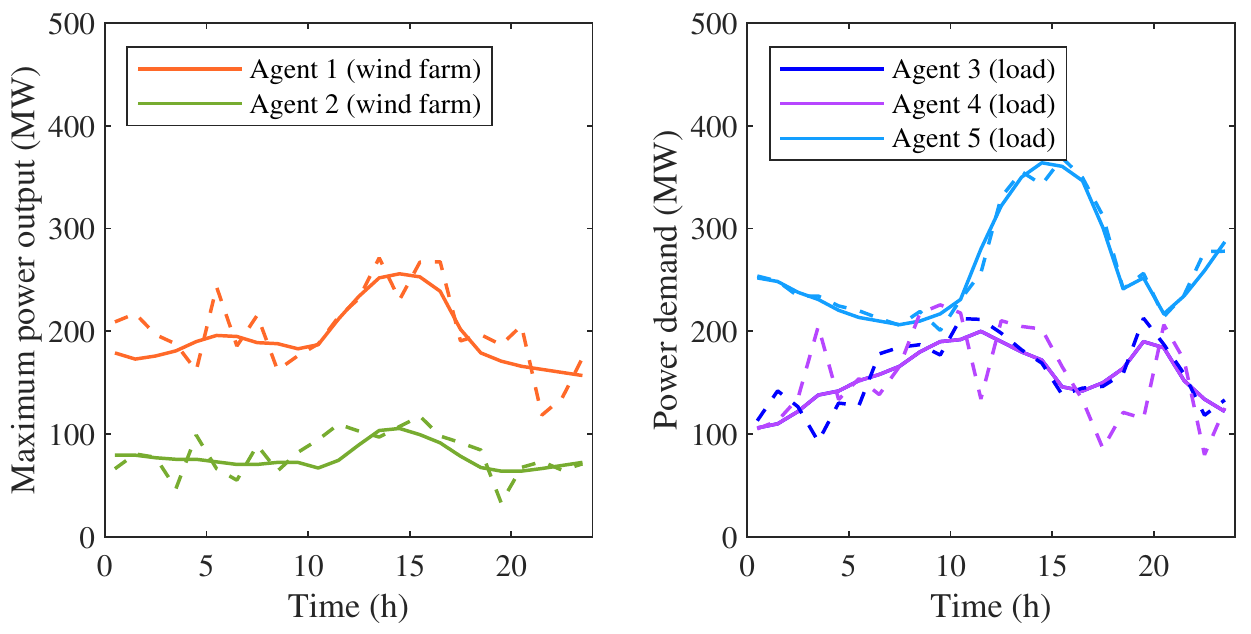}
\caption{Original forecasts (solid lines) and actual values (dashed lines).}
\label{fig:uncertainty}
\end{figure}

A 5-bus system with 3 controllable generators, 5 agents (2 wind farms and 3 uncertain loads), and 3 fixed loads is tested, whose parameters are shown in Table \ref{tab:parameter}. The time interval is 1 h. 
The operator's original forecasts $\bar u_i$ and the actual uncertainty realization $u_i$ of agent $i$ in hindsight are depicted in Fig. \ref{fig:uncertainty}. Suppose in every period $\sigma_U^2 = [8000,2000,4000,9000,1000]$ MW$^2$. 
\begin{table}[ht]
        \renewcommand{\arraystretch}{1.3}
        \renewcommand{\tabcolsep}{1em}
        \scriptsize
        \centering
        \caption{Parameters}
        \label{tab:parameter}
        \begin{tabular}{cc|cc}
        \hline
                Parameter & Value & Parameter & Value \\
               \hline
              $\rho$ & [35,30,25] \$/MWh & $T$ & 24 \\
              $P^{max}$ & [700,700,800] MW & $m$ & $1\times 10^4$ \$$\cdot$MW$^2$ \\
              $P^{min}$ & [280,280,320] MW & $\rho^c$ & 100 \$/MWh \\
              $R^\pm$ & [350,350,400] MW & $\iota$ & $10^4$ \$/MW$^2$ \\
              $\mathcal{R}^\pm$ & [350,350,400] MW & $\delta,\xi$ & $0.95,0.95$ \\
              \hline
        \end{tabular}
\end{table}

The proposed algorithm converges after 14 iterations in 987 s, which is acceptable for the day-ahead scheduling. 
The total operation cost (sum of the first-stage and the worst-case second-stage operation costs) is \$$1.115 \times 10^6$, while the prediction payments $C_i$ for the 5 agents are $\$[4.00,4.28,3.73,3.55,0] \times 10^3$, respectively. Hence, the total cost \eqref{eq:RS-UL.1} under the worst-case scenario is \$$1.131\times 10^6$. We also test the performance of the obtained day-ahead predispatch strategy when dealing with the actual uncertainty realizations. A feasible real-time redispatch strategy still exists but with a lower total cost (\$$1.071 \times 10^6$).

\begin{figure}[ht]
\centering
\includegraphics[width=1.0\columnwidth]{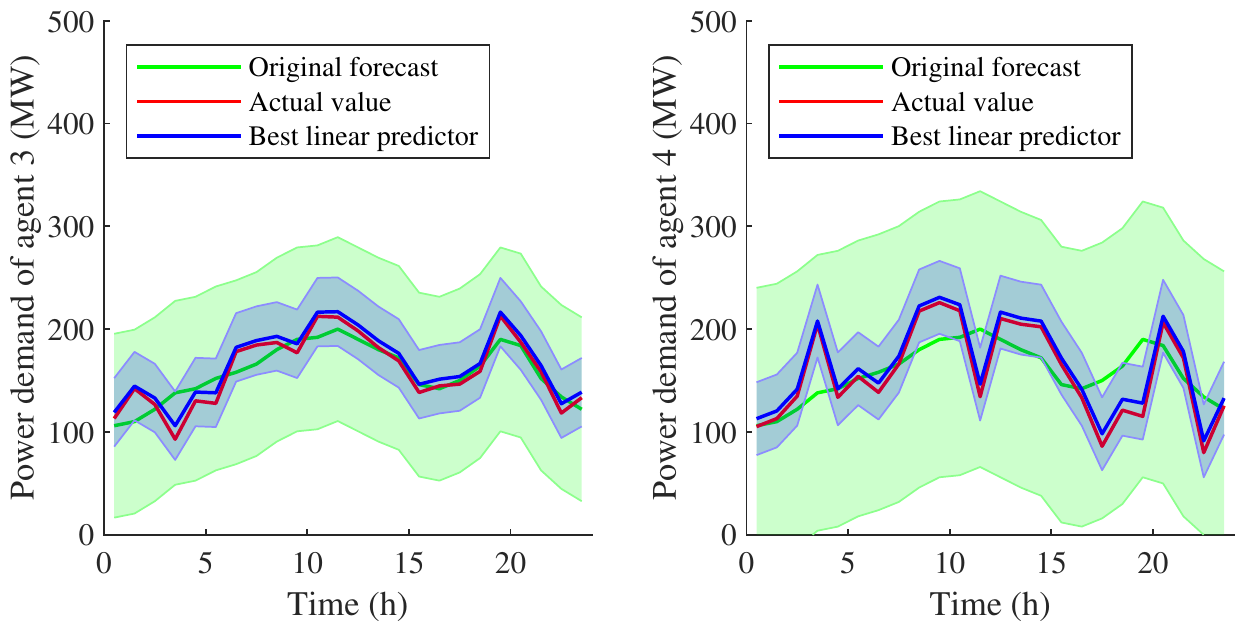}
\caption{Original and improved uncertainty sets of agents 3 and 4.}
\label{fig:region}
\end{figure}

To visualize the impact of predictions from RESs and loads on the forecast of the operator, the original (green area) and improved (blue area) uncertainty sets of agents 3 and 4 (i.e., loads 1 and 2) are shown in Fig. \ref{fig:region}. Both sets have a confidence probability $\delta = 0.95$ in each period. The centers of the original uncertainty sets are the original forecasts $\bar u_{it},\forall i,\forall t$. After purchasing the predictions from RESs and loads, the centers become the best linear predictors $u_{it}^e,\forall i,\forall t$, which are closer to the actual uncertainty realizations. The shaded areas show the variation ranges of the demand of loads $1$ and $2$.

\begin{figure}[ht]
\centering
\includegraphics[width=1\columnwidth]{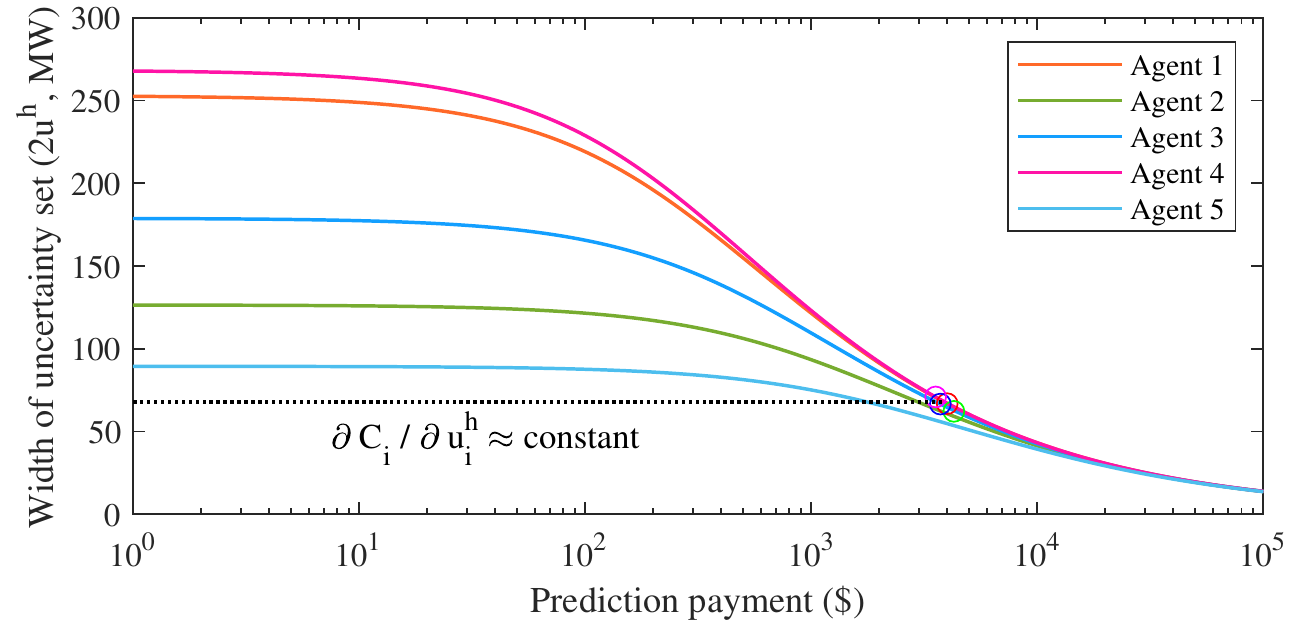}
\caption{Width of uncertainty sets under different prediction payments; circles represent the results at the optimum (the optimal payment for agent 5 is zero, which is outside of this figure).}
\label{fig:width}
\end{figure}

Both uncertainty sets contain the actual load demand, but the improved sets are much narrower so the operator is facing less uncertainty. The original forecasts of agents 3 and 4 are the same, but the uncertainty variance of agent 4 is larger, so agent 4's original uncertainty set is wider ($u_4^h =$ 134 MW) than that of agent 3 ($u_3^h =$ 89 MW). The widths of their improved sets are similar with $u_3^h=33$ MW and $u_4^h=35$ MW. Moreover, Fig. \ref{fig:width} shows how the uncertainty sets narrow as the prediction payments increase. Note that Fig. \ref{fig:width} is a semi-log plot and the prediction payment grows very quickly when the width of the uncertainty set is small. This is because the marginal prediction cost increases with a higher accuracy. In other words, it costs more to improve the accuracy of an already quite accurate prediction. The optimal payments and the corresponding widths $2 u_i^{h*},\forall i$ are also marked in Fig. \ref{fig:width}. We can find that the values of $u_i^{h*},\forall i$ are similar, which is due to the equal incremental principle, i.e., at the optimum, $\partial{C_i}/\partial{u_i^h},\forall i$ are equal. The subtle difference between $u_i^{h*},\forall i$ is caused by the linearization approximation errors. 

\subsection{Sensitivity Analysis}
We further investigate the impacts of three different factors: the agent's prediction cost coefficient $m$, the probability parameters $\delta$ and $\xi$ of the uncertainty set, and the variance of the uncertain factor $\sigma_U^2$.

\subsubsection{Impact of Prediction Cost Coefficient} 
\begin{figure}[ht]
\centering
\includegraphics[width=1.0\columnwidth]{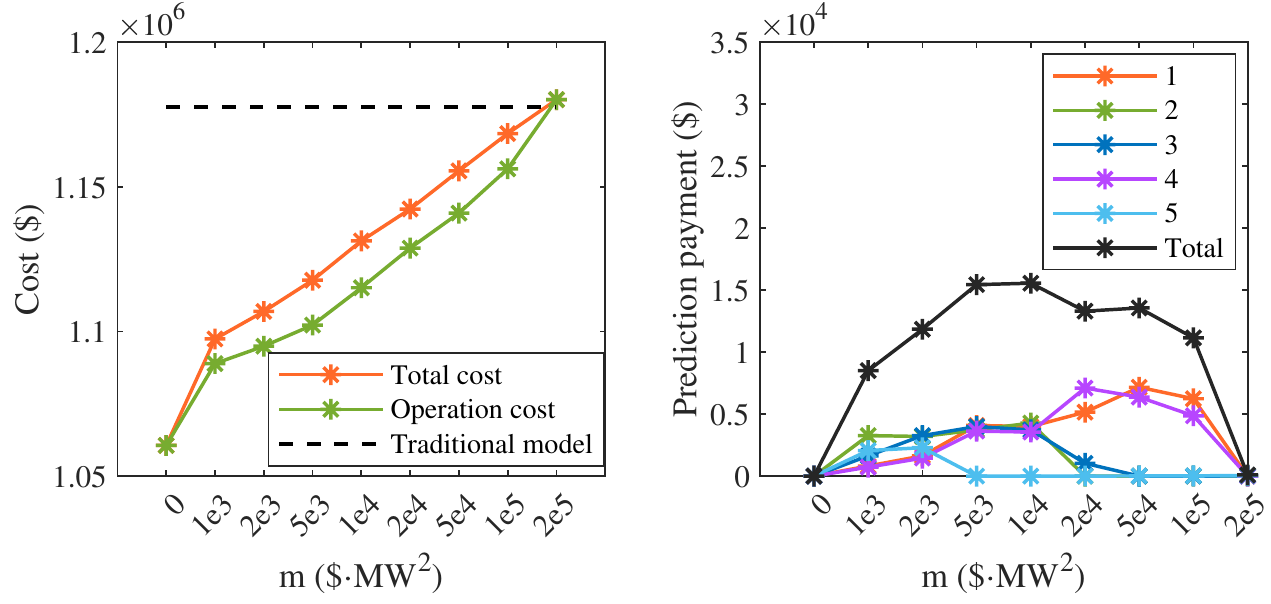}
\caption{Costs and prediction payments under different $m$.}
\label{fig:cost_payment_m}
\end{figure}

\begin{figure}[ht]
\centering
\includegraphics[width=1.0\columnwidth]{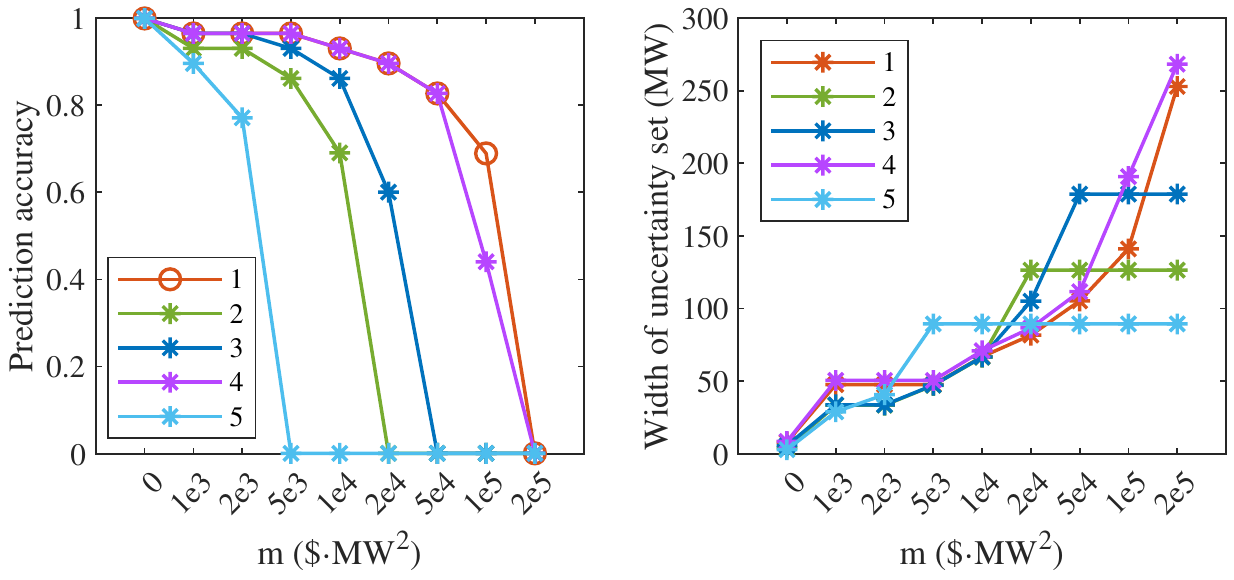}
\caption{Prediction accuracy and width of uncertainty set under different $m$.}
\label{fig:accuracy_width_m}
\end{figure}

First, we test how the strategy of the operator changes with a rising agent prediction cost by changing $m$ from 0 to $2\times 10^5$ \$$\cdot$MW$^2$. The total costs, operation costs, and prediction payments under different $m$ are shown in Fig. \ref{fig:cost_payment_m}. The change of prediction accuracy $\tau$ and the width of the improved uncertainty set are given in Fig. \ref{fig:accuracy_width_m}. We can find that when $m = 0$, the agents' prediction payments are zero since the operator can know the exact value of $u$ without making any payment, and thus, there is no uncertainty ($\tau_i=1,\forall i$). 
When $m$ is very large, e.g., $2 \times 10^5$ \$$\cdot$MW$^2$, the agents' prediction costs are extremely high, so the operator cannot afford to purchase predictions from the agents. Therefore, as shown in Fig. \ref{fig:accuracy_width_m}, the final prediction accuracy $\tau_i$ is zero for each agent $i$ and the uncertainty sets are the widest.
As $m$ grows, from Fig. \ref{fig:cost_payment_m}, the operation cost and the total cost are always less than the cost of the traditional model \eqref{eq:tr-robust} without buying predictions from the agents. The lower the $m$, the higher the operation cost reduction, showing the potential of our model.
According to Definition 2, the prediction payment $C_i$ can be interpreted as the value of prediction from agents. This value is influenced by agent's prediction cost coefficient $m$ and the system parameters. From Fig. \ref{fig:cost_payment_m}, the value of prediction of all agents follows a similar trend (first increases and then declines) and the peak value of the agent with a larger uncertainty variance $\sigma_{U_i}^2$ tends to be higher. This indicates that predictions from agents will play an increasingly important role in future power systems with higher uncertainties.


\subsubsection{Impact of Uncertainty Set Probability Parameters} 
\begin{figure}[ht]
\centering
\includegraphics[width=1.0\columnwidth]{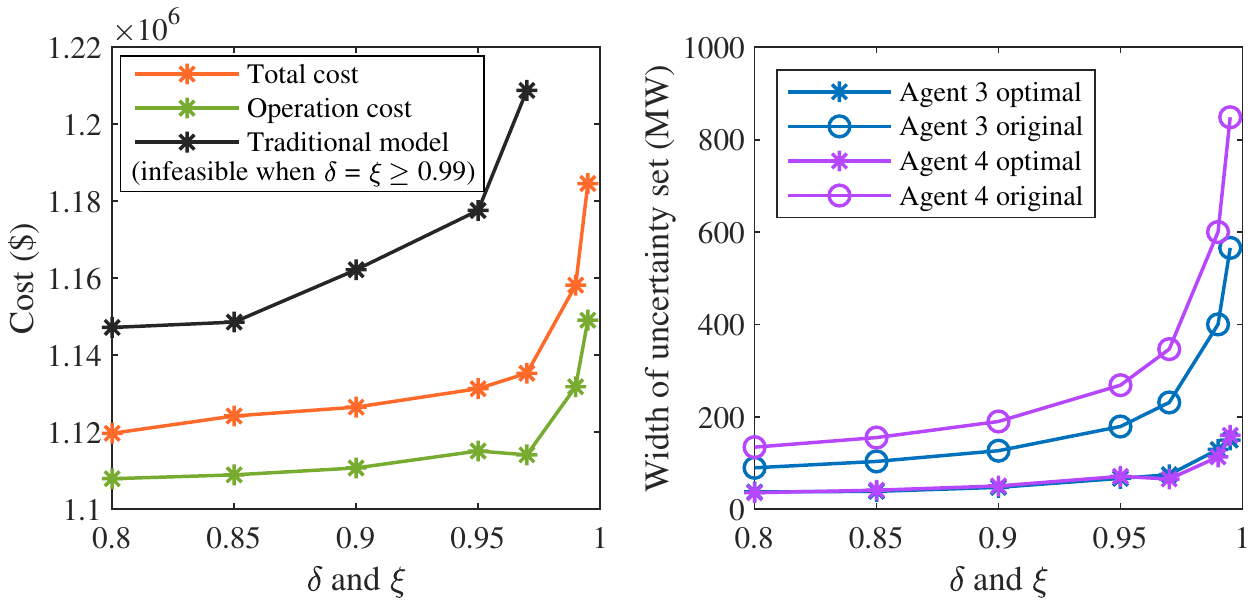}
\caption{Costs and widths of uncertainty sets under different $\delta$ and $\xi$.}
\label{fig:cost_width_xi}
\end{figure}

We next change the probability parameters $\delta$ and $\xi$ simultaneously, i.e., keeping $\delta = \xi$. The costs as well as the widths of uncertainty sets of agents 3 and 4 are shown in Fig. \ref{fig:cost_width_xi}. When $\delta$ and $\xi$ increase, the uncertainty sets expand, giving a more robust optimal predispatch strategy but also resulting in higher total and operation costs. Moreover, the total cost of the proposed model is always less than that of the traditional model without predictions from the agents (when $\delta = \xi \geq 0.99$, the traditional model become infeasible). 
The improved uncertainty sets are much narrower than the original uncertainty sets. The widths of the improved uncertainty sets of agents 3 and 4 are similar under different $\delta$ and $\xi$ due to the same reason as in Section V-A.

\subsubsection{Impact of Uncertainty Variance} 
\begin{figure}[ht]
\centering
\includegraphics[width=1.0\columnwidth]{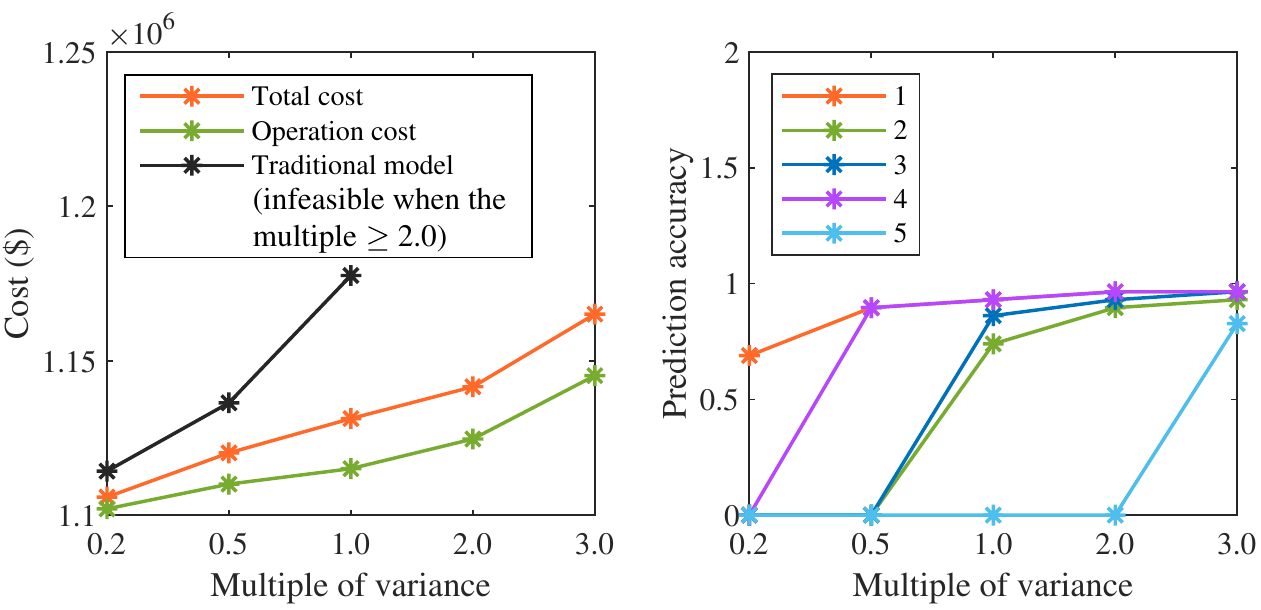}
\caption{Cost and prediction accuracy $\tau$ under different multiples of variance.}
\label{fig:cost_accuracy_sigma2}
\end{figure}

We further investigate the impact of $\sigma_U^2$, the variance of the uncertain factor (or the operator's original estimate). To do this, we multiply $\sigma_U^2$ by a positive constant. The original forecasts and actual values are still the same as those in Fig. \ref{fig:uncertainty}. The cost and prediction accuracy $\tau_i,\forall i$ are shown in Fig. \ref{fig:cost_accuracy_sigma2}. When the multiple of variance exceeds $2.0$, the traditional model is infeasible because the uncertainty is too severe. The proposed model is still feasible because buying predictions from the agents enables the operator to effectively mitigate the uncertainty they face. This shows the advantage of the proposed model.
When the variance is small, the operator already has a relatively good original estimate, so they tend to pay less for buying predictions from the agents. 
When the variance is large, at the optimum, the prediction accuracy is close to 1 for every agent, meaning that the operator relies on the predictions from the agents to make better dispatch decisions.

\subsection{Comparison With the Traditional C\&CG Algorithm}
\label{secV-C}
\begin{figure}[ht]
\centering
\includegraphics[width=1.0\columnwidth]{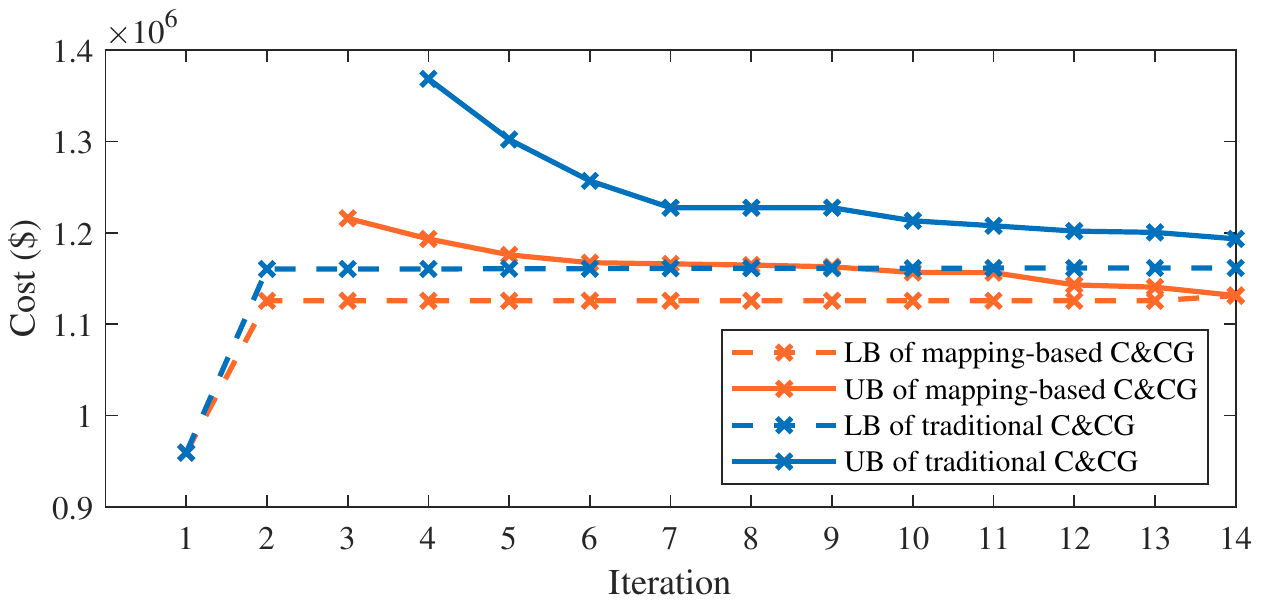}
\caption{Iteration processes of the proposed and traditional C\&CG algorithms.}
\label{fig:CCG}
\end{figure}

To show the necessity of the proposed mapping-based C\&CG algorithm, we compare it with the traditional C\&CG algorithm \cite{zeng2013solving} using the benchmark case. 
The iteration processes are shown in Fig. \ref{fig:CCG}. The mapping-based C\&CG algorithm converges to the optimal solution given in Section V-A. The traditional algorithm stops in 30 iterations with the optimal objective value equals $\$ 1.178 \times 10^6$, which is higher than that of the proposed algorithm ($\$ 1.131\times 10^6$).
This is because a previously added worst-case scenario may lie outside of the uncertainty set when the first-stage decision changes. Therefore, the master problem is no longer a relaxation of the robust optimization, which may lead to over-conservative results. Moreover, the previous scenarios that are outside of the uncertainty set hinder the improvement by buying predictions from the agents. At the optimum of the traditional C\&CG, we can find that the prediction payments are zero. Given the reasons above, the proposed algorithm is necessary.

\subsection{Out-of-Sample Test} 
To analyze the statistical performance of the obtained strategy, out-of-sample tests are conducted. To imitate various prediction errors, we randomly generate scenarios from a uniform and a Gaussian distribution with the same expectation $u$ and standard deviation, respectively. We change the standard deviation via multiplying $\sqrt{\mbox{var}(U|U^{pre})}$ by a constant from 0.5 to 2.5 and test the average total cost under the selected scenarios. Ten thousand (10000) scenarios are tested for each setting. The average total costs of the proposed algorithm and the traditional C\&CG are compared in Table \ref{tab:comparison-2}
. We can find that the proposed algorithm has lower costs. 

\begin{table}[ht]
        \renewcommand{\arraystretch}{1.3}
        \renewcommand{\tabcolsep}{1em}
        \scriptsize
        \centering
        \caption{Out-of-sample test of the proposed algorithm and the traditional C\&CG: Average total cost ($10^6$ \$)}
        \label{tab:comparison-2}
        \begin{tabular}{ccccccc}
                \hline  \multicolumn{2}{c}{Multiple of standard} & \multirow{2}{*}{0.5} & \multirow{2}{*}{1.0} & \multirow{2}{*}{1.5} & \multirow{2}{*}{2.0} & \multirow{2}{*}{2.5} \\
                \multicolumn{2}{c}{deviation} \\
               \hline
              \multirow{2}{*}{Uniform} & Proposed & 1.0733 & 1.0741 & 1.0762 & 1.0796 & 1.0840 \\
              & Traditional & 1.1237 & 1.1245 & 1.1257 & 1.1274 & 1.1294 \\
              \hline
              \multirow{2}{*}{Gaussian} & Proposed & 1.0733 & 1.0741 & 1.0762 & 1.0794 & 1.0830 \\
              & Traditional & 1.1237 & 1.1245 & 1.1257 & 1.1272 & 1.1286 \\
              \hline
        \end{tabular}
\end{table}

\subsection{Scalability}
To show the scalability of the proposed algorithm, the computational time and the number of iterations under larger systems are recorded in Table \ref{tab:comparison}. In the results, $UB/LB-1 \leq 1\%$. The time needed is all less than 2 h, which is acceptable for day-ahead generation dispatch.

\begin{table}[ht]
        \renewcommand{\arraystretch}{1.3}
        \renewcommand{\tabcolsep}{1em}
        \scriptsize
        \centering
        \caption{Computational time/number of iterations under different settings}
        \label{tab:comparison}
        \begin{tabular}{ccccc}
                \hline 
               No. of agents & 4 & 8 & 12 & 16 \\
              \hline
              33-bus (4 generators) & 84 s / 3 & 116 s / 3 & 163 s / 3 & 1576 s / 15 \\
              69-bus (8 generators) & 104 s / 3 & 329 s / 6 & 2614 s / 17 & 3609 s / 17 \\
              123-bus (8 generators) & 542 s / 6 & 968 s / 5 & 2156 s / 8 & 6411 s / 11 \\ 
              \hline
        \end{tabular}
\end{table}



\section{Conclusion}
\label{sec:conclusion}
This paper proposes a novel robust generation dispatch model in which the operator can purchase predictions from the agents to obtain a more accurate uncertainty set and make better decisions. The proposed model renders a two-stage RO with DDU. A mapping-based C\&CG algorithm with convergence guarantee is developed to solve the model. Some interesting findings are: 
\begin{itemize}
    \item Compared with the traditional model without buying predictions from the agents, the proposed model can help the operator greatly narrow the uncertainty set and reduce the total cost.
    \item The value of predictions from the agents grows with the variance of uncertainty, indicating that they will play an increasingly important role in future power systems with more volatile renewable generation.
    \item When dealing with DDU, the proposed algorithm outperforms the traditional C\&CG algorithm in terms of solution optimality.
\end{itemize}
A detailed predictive information market design between the operator and agents and a more efficient solution algorithm will be our future research directions.



\ifCLASSOPTIONcaptionsoff
\newpage
\fi

\bibliographystyle{IEEEtran}
\bibliography{IEEEabrv,mybib}

\appendices
\makeatletter
\@addtoreset{equation}{section}
\@addtoreset{theorem}{section}
\makeatother
\setcounter{equation}{0}  
\renewcommand{\theequation}{A.\arabic{equation}}
\renewcommand{\thetheorem}{A.\arabic{theorem}}
\section{Proof of Lemma \ref{lemma-1}}
\label{apendix-A}
Set a function
\begin{align}
    g(\alpha_i, \beta_i) :=~ & \mathbbm{E} \left[ \left( U_i - (\alpha_i + \beta_i U_i^{pre}) \right)^2 \right] \nonumber \\
    =~ & \alpha_i^2 + \mathbbm{E}[(U_i^{pre})^2] \beta_i^2 + 2 \mathbbm{E}[U_i^{pre}] \alpha_i \beta_i \nonumber \\
    & - 2 \mathbbm{E}[U_i] \alpha_i - 2 \mathbbm{E}[U_i U_i^{pre}] \beta_i + \mathbbm{E}[U_i^2].
\end{align}

Note that $\epsilon_i$ is independent of $U_i$, then
\begin{align}
    & \mathbbm{E}[U_i^{pre}] = \mathbbm{E}[U_i - \epsilon_i] = \mathbbm{E}[U_i] = \overline{u}_i, \nonumber \\
    & \mathbbm{E}[U_i^2] = \mbox{var}[U_i] + (\mathbbm{E}[U_i])^2 = \sigma_{U_i}^2 + \overline{u}_i^2, \nonumber \\
    & \mathbbm{E}[U_i \epsilon_i] = \mathbbm{E}[U_i] \mathbbm{E}[\epsilon_i] = 0, \nonumber \\
    & \mathbbm{E}[U_i U_i^{pre}] = \mathbbm{E}[U_i^2] - \mathbbm{E}[U_i \epsilon_i] = \mathbbm{E}[U_i^2] = \sigma_{U_i}^2 + \overline{u}_i^2, \nonumber \\
    & \mathbbm{E}[(U_i^{pre})^2] = \mathbbm{E}[U_i^2] - 2 \mathbbm{E}[U_i \epsilon_i] + \mathbbm{E}[\epsilon_i^2] = \sigma_{U_i}^2 + \overline{u}_i^2 + \sigma_{\epsilon_i}^2. \nonumber
\end{align}
Since $\alpha_i$ and $\beta_i$ minimizes $g(\alpha_i, \beta_i)$, we have
\begin{align}
    \frac{\partial g(\alpha_i, \beta_i)}{\partial \alpha_i} = 0,
    \frac{\partial g(\alpha_i, \beta_i)}{\partial \beta_i} = 0,
\end{align}
whose solution is 
\begin{align}
    \beta_i=\frac{\sigma^2_{U_i}}{\sigma^2_{U_i}+\sigma^2_{\epsilon_i}}, ~\alpha_i=(1-\beta_i)\overline{u}_i.
\end{align}
Then
\begin{align}
    \mathbbm{E}[U_i-U_i^e] = \overline{u}_i - \alpha_i - \beta_i \overline{u}_i = 0, \nonumber
\end{align}
and
\begin{align}
    & \mathrm{cov}(U_i-U_i^e,U_i^{pre}) \nonumber \\
    =~ & \mathbbm{E}[(U_i-U_i^e)U_i^{pre}] - \mathbbm{E}[U_i-U_i^e] \mathbbm{E}[U_i^{pre}] \nonumber \\
    =~ & \mathbbm{E}[(U_i-(\alpha_i + \beta_i U_i^{pre}))U_i^{pre}] \nonumber \nonumber \\
    =~ & \mathbbm{E}[U_i U_i^{pre}] - \alpha_i \mathbbm{E}[U_i^{pre}] - \beta_i \mathbbm{E}[(U_i^{pre})^2] = 0. \nonumber
\end{align}
This completes the proof.

\setcounter{equation}{0}  
\renewcommand{\theequation}{B.\arabic{equation}}
\renewcommand{\thetheorem}{B.\arabic{theorem}}
\section{Proof of Lemma \ref{lemma-2}}
\label{apendix-B}

Since $U_i^e = \alpha_i + \beta_i U_i^{pre}$ is a function of $U_i^{pre}$, for any random variable $X$ with finite second moment we have $\mathbbm{E}[X U_i^e | U_i^{pre}] = U_i^e \mathbbm{E}[X|U_i^{pre}]$ \cite[p. 348]{grimmett2001probability}. Therefore, 
\begin{align}
    & \mbox{var}[\eta_i|U_i^{pre}] \nonumber \\
    =~ & \mathbbm{E}[(U_i - U_i^e)^2|U_i^{pre}] - (\mathbbm{E}[U_i - U_i^e|U_i^{pre}])^2 \nonumber \\
    =~ & (\mathbbm{E}[U_i^2|U_i^{pre}] - 2 U_i^e \mathbbm{E}[U_i|U_i^{pre}] + (U_i^e)^2 \nonumber) \\
    & - ((\mathbbm{E}[U_i|U_i^{pre}])^2 - 2 U_i^e \mathbbm{E}[U_i|U_i^{pre}] + (U_i^e)^2) \nonumber \\
    =~ & \mathbbm{E}[U_i^2|U_i^{pre}] - (\mathbbm{E}[U_i|U_i^{pre}])^2 \nonumber \\
    =~ & \mbox{var}[U_i|U_i^{pre}].
\end{align}

Moreover, we have $\eta_i$ and $U_i^{pre}$ are independent random variables. Therefore,
\begin{align}
    \mbox{var}[U_i|U_i^{pre}]=~ & \mbox{var}[\eta_i|U_i^{pre}] = \mbox{var}[\eta_i] \nonumber\\
    = ~ & \mbox{var}[U_i-(1-\beta_i)\overline{u}_i-\beta_i U_i^{pre}] \nonumber\\
    = ~ & \mbox{var}[(1-\beta_i)U_i+\beta_i\epsilon_i]\nonumber\\
    \label{eq:lemma-2-eta}
    = ~ & (1-\beta_i)^2 \sigma_{U_i}^2 + \beta_i^2 \sigma_{\epsilon_i}^2.
\end{align}

Moreover,
\begin{align}
    \mbox{var}[U_i|U_i^{pre}] = \sigma_{U_i}^2 \frac{\sigma_{\epsilon_i}^2}{\sigma_{\epsilon_i}^2+\sigma_{U_i}^2} \le \sigma_{U_i}^2.
\end{align}
This completes the proof.

\setcounter{equation}{0}  
\renewcommand{\theequation}{C.\arabic{equation}}
\renewcommand{\thetheorem}{C.\arabic{theorem}}
\section{Proof and Discussion of Lemma \ref{lemma-3}}
\label{apendix-C}
\subsection{Proof of Lemma \ref{lemma-3}}
Recall that according to Lemma \ref{lemma-1} we have $\mathbbm{E}[\eta_{it}] = \mathbbm{E}[U_{it}-U_{it}^e] = 0$. Combining the fact $\mbox{var}[\eta_{it}]=\mbox{var}[U_{it}|U_{it}^{pre}]$ shown in \eqref{eq:lemma-2-eta}, we have $\mathbbm{E}[\eta_{it}/u_{it}^h] = 0$ and $\mbox{var}[\eta_{it}/u_{it}^h] = 1 - \delta$. Then according to the Chebyshev inequality,
\begin{align}
    ~ & \mathcal{P}\left(|\eta_{it}/u_{it}^h|\ge 1\right) \nonumber\\
    = ~ & \mathcal{P}\left(|(\eta_{it}/u_{it}^h) - \mathbbm{E}[\eta_{it}/u_{it}^h]|\ge \sqrt{\mbox{var}[\eta_{it}/u_{it}^h]/(1-\delta)}\right) \nonumber\\
    \le ~ & 1/\left(\sqrt{1/(1-\delta)}\right)^2 = 1-\delta.
\end{align}
Note that $v_{it} = |\eta_{it}/u_{it}^h|$. Thus, we have $\mathcal{P}(v_{it}\ge 1)=\mathcal{P}(|\eta_{it}/u_{it}^h|\ge 1)\le 1-\delta$.

Again by the Chebyshev inequality,
\begin{align}
    ~ & \mathcal{P}\left(\sum \nolimits_i v_{it}\ge\mathbbm{E}[\sum \nolimits_i v_{it}] + \sqrt{\mbox{var}[\sum \nolimits_i v_{it}]/(1-\xi)} \right) \nonumber\\
    = ~ & \mathcal{P}\left(\sum \nolimits_i v_{it}-\mathbbm{E}[\sum \nolimits_i v_{it}] \ge \sqrt{\mbox{var}[\sum \nolimits_i v_{it}]}\sqrt{1/(1-\xi)}\right) \nonumber\\
    \le ~ & \left\{1/(1-\xi)\right\}^{-1} = 1-\xi .
\end{align}

We find an upper bound for $\mathbbm{E}[\sum \nolimits_i v_{it}] + \sqrt{\mbox{var}[\sum \nolimits_i v_{it}]/(1-\xi)}$. $\{\eta_{it}/\mbox{var}[\eta_{it}],\forall i,\forall t\}$ are independent and identically distributed (i.i.d.), so $\{v_{it},\forall i,\forall t\}$ are also i.i.d. random variables. Therefore, $\mathbbm{E}[\sum_i v_{it}]=I \mathbbm{E}[v_{it}]$ and $\mbox{var}[\sum_i v_{it}]=I \mbox{var}[v_{it}]$. Moreover,
\begin{align}
    \mbox{var}[v_{it}] =~ & \mathbbm{E}[v_{it}^2] - (\mathbbm{E}[v_{it}])^2 \nonumber \\
    =~ & \mathbbm{E}[(\eta_{it}/u_{it}^h)^2] - (\mathbbm{E}[v_{it}])^2 \nonumber \\
    =~ & \mbox{var}[\eta_{it}/u_{it}^h] + (\mathbbm{E}[\eta_{it}/u_{it}^h])^2 - (\mathbbm{E}[v_{it}])^2 \nonumber \\
    =~ & 1-\delta - (\mathbbm{E}[v_{it}])^2.
\end{align}
Then
\begin{align}
    &  \mathbbm{E}[\sum \nolimits_i v_{it}] + \sqrt{\mbox{var}[\sum \nolimits_i v_{it}]/(1-\xi)} \nonumber \\
    =~ & I \mathbbm{E}[v_{it}] + \sqrt{I(1-\delta - (\mathbbm{E}[v_{it}])^2)/(1-\xi)} ~=: G(\mathbbm{E}[v_{it}]) \nonumber
\end{align}
is a function of $\mathbbm{E}[v_{it}]$ for $0 \le \mathbbm{E}[v_{it}] \le \sqrt{1-\delta}$. By calculating its derivative, it is easy to show that $G(\mathbbm{E}[v_{it}])$ first increases and then declines, whose unique maximum value is $\Gamma_S$.
Therefore, we have $\mathcal{P}(\sum_i v_{it} \ge \Gamma_S) \le 1-\xi,\forall t$. Similarly, we can prove that $\mathcal{P}(\sum_t v_{it} \ge \Gamma_T) \le 1-\xi,\forall i$.

\subsection{Discussion of Lemma \ref{lemma-3}}

The bound $\mathcal{P}(v_{it} \geq 1) \leq 1-\delta$ is tight for $\eta_{it}/\sqrt{\mbox{var}[\eta_{it}]}$ with the following discrete probability distribution:
\begin{align}
    & \mathcal{P}(\eta_{it}/\sqrt{\mbox{var}[\eta_{it}]}=-1/\sqrt{1-\delta}) = (1-\delta)/2, \nonumber \\
    & \mathcal{P}(\eta_{it}/\sqrt{\mbox{var}[\eta_{it}]}=0)=\delta, \nonumber \\
    & \mathcal{P}(\eta_{it}/\sqrt{\mbox{var}[\eta_{it}]}=1/\sqrt{1-\delta}) = (1-\delta)/2, \nonumber 
\end{align}
which satisfies 
\begin{align}
    \mathbb{E}[\eta_{it}/\sqrt{\mbox{var}[\eta_{it}]}] = 0,~ \mbox{var}[\eta_{it}/\sqrt{\mbox{var}[\eta_{it}]}] = 1. \nonumber
\end{align}

Now we focus on the bound $\mathcal{P}(\sum_i v_{it} \geq \Gamma_S) \leq 1-\xi$, where $\Gamma_S$ is given by \eqref{eq:uncertainty-set-3}. We consider some specific testing probability distributions of $\eta_{it}/\sqrt{\mbox{var}[\eta_{it}]}$ and calculate $\mathcal{P}(\sum_i v_{it} \geq \Gamma_S)$ under varying $\Gamma_S$. 

First, consider the probability distribution with
\begin{align}
    \mathcal{P}(\eta_{it}/\sqrt{\mbox{var}[\eta_{it}]}=1) =  \mathcal{P}(\eta_{it}/\sqrt{\mbox{var}[\eta_{it}]}=-1) = 1/2. \nonumber
\end{align}
Then $\mathcal{P}(v_{it}=\sqrt{1-\delta}) = 1, \forall i$ and $\mathcal{P}(\sum_i v_{it} = I \sqrt{1-\delta}) = 1$. Therefore, $\mathcal{P}(\sum_i v_{it} \geq \Gamma_S) = 1$ for $\Gamma_S \leq I\sqrt{1-\delta}$ and $\mathcal{P}(\sum_i v_{it} \geq \Gamma_S) = 0$ for $\Gamma_S > I\sqrt{1-\delta}$.

Second, construct a discrete probability distribution as follows with parameter $\xi \in (0,1)$.
\begin{align}
    & \mathcal{P}(\eta_{it}/\sqrt{\mbox{var}[\eta_{it}]}=a)=\mathcal{P}(\eta_{it}/\sqrt{\mbox{var}[\eta_{it}]}=-a)=p, \nonumber \\
    & \mathcal{P}(\eta_{it}/\sqrt{\mbox{var}[\eta_{it}]}=0)=1-2p, \nonumber
\end{align}
where
\begin{align}
    a = \sqrt{\frac{1+I-I\xi}{I(1-\xi)}},~ p = \frac{I(1-\xi)}{2(1+I-I\xi)}. \nonumber
\end{align}
Then $\mathcal{P}(v_{it}=a\sqrt{1-\delta}) = 2p$ and $\mathcal{P}(v_{it}=0)=1-2p,\forall i$, based on which $\mathcal{P}(\sum_i v_{it} \geq \Gamma_S)$ can be calculated under different $\Gamma_S$.

\begin{figure}[ht]
\centering
\includegraphics[width=1.0\columnwidth]{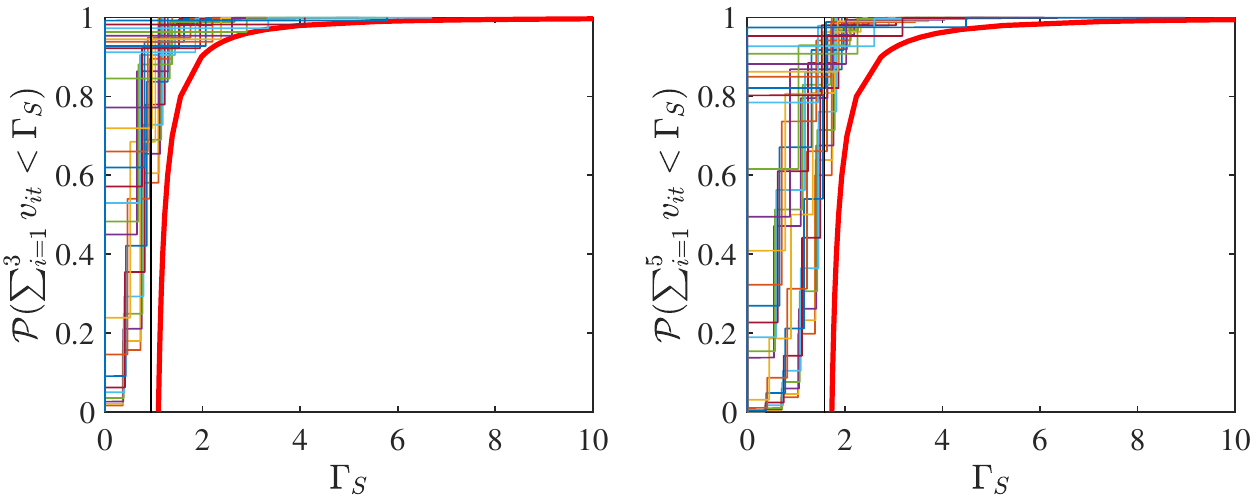}
\caption{Bounds \eqref{eq:uncertainty-set-3} (red thick curves) and testing values (other curves) of $\mathcal{P}(\sum_i v_{it} \geq \Gamma_S)$ when $I=3$ (left) and $I=5$ (right).}
\label{fig:lemma3}
\end{figure}

Bound \eqref{eq:uncertainty-set-3} is compared with the values of $\mathcal{P}(\sum_i v_{it} \geq \Gamma_S)$ under the two kinds of testing probability distributions in Fig. \ref{fig:lemma3}, where $I = 3$ or $I = 5$ and $\xi$ varies in $(0,1)$. All the curves of the testing probability distributions are above the curves of bound \eqref{eq:uncertainty-set-3}, which verifies the effectiveness of the bound.

When $I$ is large, by the central limit theorem, the probability distribution of $\sum_i v_{it}$ can be approximated by a normal distribution with mean $I\mathbb{E}[v_{it}]$ and variance $I\mbox{var}[v_{it}] = I(1-\delta-(\mathbb{E}[v_{it}])^2)$. Therefore, with a high probability, $\sum_i v_{it}$ will not exceed $I\sqrt{1-\delta}$ much. Hence, $\Gamma_S \approx I\sqrt{1-\delta}$ can be used as the bound in $\mathcal{P}(\sum_i v_{it} \geq \Gamma_S) \approx 0$ when $I$ is large.

\setcounter{equation}{0}  
\renewcommand{\theequation}{D.\arabic{equation}}
\renewcommand{\thetheorem}{D.\arabic{theorem}}
\section{Proof of Theorem \ref{thm-1}}
\label{apendix-D}

Suppose the optimal solution of the two-stage RO model \eqref{eq:new-robust} is $(x^*,\tau^*)$ and the optimal objective value is
\begin{align}
    O^* := \min_{x \in \mathcal{X} \cap \tilde{\mathcal{X}}_R \atop \tau_i \in [0,1], \eqref{eq:precost}, \forall i \in \mathcal{I}} f(x) + \sum_{i \in \mathcal{I}} C_i(\tau_i) + \max_{u \in \mathcal{U}(\tau)} \min_{y \in \mathcal{Y}(x,u)} g(y). \nonumber
\end{align}

We start the proof of Theorem \ref{thm-1} by providing the following claims. For any $K \in \mathbb{Z}^{+}$:

(1) \textbf{Claim 1}: $LB_K \le O^* \le UB_K$;

(2) \textbf{Claim 2}: If the algorithm does not terminate after $K$ iterations, then for any $K_1, K_2 \in [K]$, we have $\phi_{K_1}^* \ne \phi_{K_2}^*$.
\vspace{0.5em}

\emph{Proofs of claims}:

(1) \textbf{Claim 1}: The master problem in the $K$-th iteration is equivalent to
\begin{align}
    {LB}_{K} = \!\!\!\!\!\min_{x \in \mathcal{X} \cap \mathcal{X}_{K-1} \atop \tau_i \in [0,1], \eqref{eq:precost}, \forall i \in \mathcal{I}} \!\!\!\!\! f(x) + \!\!\sum_{i \in \mathcal{I}} C_i(\tau_i) + \!\!\!\max_{u \in \mathcal{U}_{K-1}(\tau)} \min_{y \in \mathcal{Y}(x,u)} g(y), \nonumber
\end{align}
where
\begin{align}
    \mathcal{U}_{K-1}(\tau) :=~ & \{ u^k = u^e(\tau) + u^h(\tau) \phi_k^* ~|~ k \in [K-1] \}, \\
    \mathcal{X}_{K-1}:=~ & \{ x~|~\mathcal{Y}(x,u) \neq \emptyset, \forall u \in \mathcal{U}_{K-1}(\tau) \}.
\end{align}
Since $\mathcal{U}_{k-1}(\tau)  \subset  \mathcal{U}(\tau)$ and $\mathcal{X}_{K-1}  \supset \tilde{\mathcal{X}}_R$, we have $LB_K \le O^*$.

Next, we prove $UB_K \ge O^*$ by induction. First of all, $UB_0=+\infty \ge O^*$. Suppose for the sake of induction that $UB_{K-1} \ge O^*$, then if $1^{\top} s_K^*>0$, we have $UB_K=UB_{K-1} \ge O^*$; otherwise, $(x_{K}^*, \tau_K^*)$ is robust feasible and
\begin{align}
    {UB}_{K} =~ & f(x_{K}^*) + \sum_i C_{i,K}^* + g(y_{K}^*) \nonumber \\
    =~ & f(x_{K}^*) + \sum_i C_{i,K}^* + \max_{u \in \mathcal{U}(\tau_{K}^*)} \min_{y \in \mathcal{Y}(x_{K}^*,u)} g(y) \nonumber \\
    \ge~ & O^*. \nonumber
\end{align}
The last inequality is due to the optimality of $(x^*,\tau^*)$.

(2) \textbf{Claim 2}: Without loss of generality, we assume that $K_1<K_2$. If we have $\phi_{K_1}^*=\phi_{K_2}^*$, then $u_{K_2}^* = u^e(\tau_{K_2}^*)+u^h(\tau_{K_2}^*)\phi_{K_2}^* =u^e(\tau_{K_2}^*)+u^h(\tau_{K_2}^*)\phi_{K_1}^* \in \mathcal{U}_{K_2-1}(\tau_{K_2}^*)$, so $x_{K_2}^*$ must be robust feasible. Moreover,
\begin{align}
    LB_{K_2} =~ & f(x_{K_2}^*)+\sum_i C_{i,K_2}^* + \max_{u \in \mathcal{U}_{K_2-1}(\tau_{K_2}^*)} \min_{y \in \mathcal{Y}(x_{K_2}^*, u)} g(y) \nonumber\\
    \ge ~ & f(x_{K_2}^*)+\sum_i C_{i,K_2}^* + \min_{y \in \mathcal{Y}(x_{K_2}^*, u_{K_2}^*)} g(y) = UB_{K_2}.
\end{align}

Together with $LB_{K_2} \le UB_{K_2}$ from Claim 1, we have $LB_{K_2}=UB_{K_2}$. This contradicts to the assumption that the algorithm does not terminate after $K\ge K_2$ iterations.

\vspace{0.5em}

Now the proof of \textbf{Theorem 1} is given below.

First, we prove that the algorithm converges in $\mathcal{O}(n_U)$ iterations. With Lemma \ref{lemma-4}, we know that the worst-case scenario $u_k^*$ can be achieved at a vertex of $\mathcal{U}_k(\tau_k^*)$. A vertex of the set $\Phi$ corresponds to a vertex of the set $\mathcal{U}_k(\tau_k^*)$. According to Claim 2, the same vertex of $\Phi$ will not be picked up twice. Moreover, the number of vertices of $\Phi$ is $n_U$. Hence, the algorithm stops in $\mathcal{O}(n_U)$ iterations.

Suppose the algorithm terminates after $K \le n_U$ iterations.

Next, we show the robust feasibility of $(x_{K}^*,\tau_{K}^*)$. Obviously, we have $x_K^* \in \mathcal{X}$ and $\tau_K^* \in [0,1]$. Moreover, $s_K^*=0$ when the algorithm terminates, so $(x_{K}^*,\tau_{K}^*)$ is robust feasible.

Finally, we show the optimality of $(x_K^*,\tau_K^*)$. According to Claim 1, we have $LB_K \le O^* \le UB_K$. Together with the condition for termination $|UB_K-LB_K| \le \epsilon$, we have 
\begin{align}
    |UB_{K}-O^*| \le |UB_K-LB_K| \le \epsilon.
\end{align}

This completes the proof.

\end{document}